\documentclass[12pt]{amsart}
\usepackage[colorlinks]{hyperref}
\usepackage{enumerate}
\input{diagrams.tex}
\def\pdfsyncstart{}
\def\pdfsyncstop{}

\def\bdi{\pdfsyncstop\begin{diagram}}
\def\edi{\end{diagram}\pdfsyncstart}

\theoremstyle{plain}

\newtheorem{thm}{Theorem}[section]
\newtheorem{cor}[thm]{Corollary}
\newtheorem{lem}[thm]{Lemma}
\newtheorem{prop}[thm]{Proposition}

\theoremstyle{definition}
\newtheorem{defi}[thm]{Definition}
\newtheorem{defis}[thm]{Definitions}
\newtheorem{conj}[thm]{Conjecture}
\newtheorem{conv}[thm]{Convention}
\newtheorem{nota}[thm]{Notation}
\newtheorem{rem}[thm]{Remark}
\newtheorem{rems}[thm]{Remarks}
\newtheorem{exa}[thm]{Example}
\newtheorem{exas}[thm]{Examples}
\newtheorem{sit}[thm]{}

\newcommand{\brem}{\begin{rem}}
\newcommand{\brems}{\begin{rems}}
\newcommand{\erem}{\end{rem}}
\newcommand{\erems}{\end{rems}}
\newcommand{\bexa}{\begin{exa}}
\newcommand{\bexas}{\begin{exas}}
\newcommand{\eexa}{\end{exa}}
\newcommand{\eexas}{\end{exas}}
\newcommand{\bdefi}{\begin{defi}}
\newcommand{\edefi}{\end{defi}}
\newcommand{\bdefis}{\begin{defis}}
\newcommand{\edefis}{\end{defis}}
\newcommand{\bcor}{\begin{cor}}
\newcommand{\ecor}{\end{cor}}
\newcommand{\blem}{\begin{lem}}
\newcommand{\elem}{\end{lem}}
\newcommand{\bconv}{\begin{conv}}
\newcommand{\econv}{\end{conv}}
\newcommand{\bconj}{\begin{conj}}
\newcommand{\econj}{\end{conj}}
\newcommand{\bprop}{\begin{prop}}
\newcommand{\eprop}{\end{prop}}
\newcommand{\bthm}{\begin{thm}}
\newcommand{\ethm}{\end{thm}}
\newcommand{\bnota}{\begin{nota}}
\newcommand{\enota}{\end{nota}}
\newcommand{\bsit}{\begin{sit}}
\newcommand{\esit}{\end{sit}}
\newcommand{\be}{\begin{equation}}
\newcommand{\ee}{\end{equation}}
\newcommand{\bproof}{\begin{proof}}
\newcommand{\eproof}{\end{proof}}
\def\ba{\begin{array}}
\def\ea{\end{array}}

\def\dashto{\dashrightarrow}
\def\lto{\longrightarrow}

\def\st{{\rm st}}
\def\reg{{\rm reg}}
\def\ext{{\rm ext}}
\def\no{\noindent}

\newcommand{\Spec}{\operatorname{Spec}}
\newcommand{\Div}{\operatorname{Div}}
\newcommand{\Pic}{\operatorname{Pic}}

\newcommand{\Frac}{\operatorname{Frac}}
\newcommand{\Proj}{\operatorname{Proj}}

\newcommand{\Aut}{{\operatorname{Aut}}}
\newcommand{\ML}{{\operatorname{ML}}}

\newcommand{\supp}{{\operatorname{supp}}}
\renewcommand{\div}{{\operatorname{div}}}

\newcommand{\SL}{{\bf {SL}}}
\newcommand{\PGL}{{\bf {PGL}}}

\def\fF{{\mathfrak F}}

\def\cO{{\mathcal O}}

\def\bC{{\bar C}}
\def\bD{{\bar D}}
\def\bV{{\bar V}}

\def\tC{{\tilde C}}
\def\tO{{\tilde O}}
\def\tD{{\tilde D}}

\def\tV{{\tilde V}}

\def\hC{{\hat C}}

\def\tgamma{{\tilde\gamma}}

\newcommand{\la}{\label}

\newcommand{\A}{{\mathbb A}}
\newcommand{\PP}{{\mathbb P}}

\newcommand{\C}{{\mathbb C}}
\newcommand{\Q}{{\mathbb Q}}
\newcommand{\Z}{{\mathbb Z}}
\newcommand{\N}{{\mathbb N}}
\newcommand{\T}{{\mathbb T}}

\newcommand{\G}{{\Gamma}}

\newcommand{\p}{{\partial}}

\def\bO{{\bar O}}

\newcommand{\nlin}{\unitlength1mm\begin{picture}(0,9.25)
                      \put(0,0.75){\line(0,1){8.5}}
                     \end{picture}}

\newcommand{\vlin}[1]{\hspace{0.75mm}\unitlength1mm\begin{picture}
(#1,0)
                      \put(0,0){\line(1,0){#1}}
                     \end{picture}\hspace{0.75mm}\rule[-3mm]{0mm}
                     {4mm}}
\newcommand{\vpoi}{\unitlength1mm\begin{picture}(0,10)
                      \multiput(0,3)(0,2){3}{\circle*{0.15}}
                    \end{picture}}

\def\llin{\vlin{11.5}}
\newcommand{\lin}{\vlin{8.5}}

\newcommand{\co}[1]{\unitlength1mm\begin{picture}(0,8)
   \put(0,0){\circle{1.5}}
   \put(0,3){\makebox(0,5)[b]{$#1$}}
                     \end{picture}}
\newcommand{\mybox}{\unitlength1mm\begin{picture}(0,1.5)
   \put(-0.75,-0.75){\line(0,1){1.5}}
   \put(-0.75,-0.75){\line(1,0){1.5}}
   \put(0.75,0.75){\line(0,-1){1.5}}
   \put(0.75,0.75){\line(-1,0){1.5}}
   \end{picture}}
\newcommand{\boxo}[1]{\unitlength1mm\begin{picture}(0,8)
   \put(0,0){\mybox}
   \put(0,3){\makebox(0,5)[b]{$#1$}}
                     \end{picture}}
\newcommand{\xbox}{\unitlength1mm\begin{picture}(0,1.5)
   \put(0,0){$\mybox$}
   \put(-0.75,0){\line(1,0){1.5}}
   \put(0,-0.75){\line(0,1){1.5}}
   \end{picture}}
\newcommand{\xboxo}[1]{\unitlength1mm\begin{picture}(0,8)
   \put(0,0){\xbox}
   \put(0,3){\makebox(0,5)[b]{$#1$}}
                     \end{picture}}
\newcommand{\cou}[2]{\unitlength1mm\begin{picture}(0,8)
   \put(0,0){\circle{1.5}}
   \put(0,3){\makebox(0,5)[b]{$#1$}}
   \put(0,-7){\makebox(0,4)[t]{$#2$}}
     \end{picture}
     \rule[-7mm]{0mm}{7mm}}

\newcommand{\crl}[2]{\unitlength1mm\begin{picture}(0,8)
   \put(0,0){\circle{1.5}}
   \put(-5,0){\makebox(0,5)[b]{$#1$}}
  \put(5,0){\makebox(0,5)[b]{$#2$}}
     \end{picture}
     \rule[-7mm]{0mm}{7mm}}

\newcommand{\cshiftup}[2]{\unitlength1mm\begin{picture}(0,9.25)
                      \put(0,10){\crl{#1}{#2}}
                     \end{picture}}

\newcommand{\xbrl}[2]{\unitlength1mm\begin{picture}(0,8)
   \put(0,0){\xbox}
   \put(-5,0){\makebox(0,5)[b]{$#1$}}
  \put(5,0){\makebox(0,5)[b]{$#2$}}
     \end{picture}
     \rule[-7mm]{0mm}{7mm}}

\newcommand{\xbshiftup}[2]{\unitlength1mm\begin{picture}(0,9.25)
                      \put(0,10){\xbrl{#1}{#2}}
                     \end{picture}}

\newcommand{\rbh}[1]{\raisebox{1mm}[1mm][-1mm]{#1}}

\title{Completions of $\C^*$-surfaces}

\author{Hubert Flenner}
\address{Fakult\"at f\"ur Mathematik,
Ruhr Universit\"at Bochum, Geb.\ NA 2/72, Universit\"ats\-str.\
150, 44780 Bochum, Germany}
\email{Hubert.Flenner@rub.de}

\author{Shulim Kaliman}
\address{Department of Mathematics,
University of Miami, Coral Gables, FL  33124, U.S.A.}
\email{kaliman@math.miami.edu}

\author{Mikhail Zaidenberg}
\address{Universit\'e
Grenoble I, Institut Fourier, UMR 5582 CNRS-UJF, BP 74, 38402
St.\ Martin d'H\`eres c\'edex, France}
\email{zaidenbe@ujf-grenoble.fr}

\thanks{
{\bf Acknowledgements:} This research was started during a visit
of the first two authors at the Fourier Institute, Grenoble, and
continued during a visit of two of us at the Max Planck Institute
of Mathematics in Bonn. We thank these institutions for their
generous support and excellent working conditions. The research
was also partially supported by NSA Grant no. MDA904-00-1-0016 and
the DFG-Schwerpunkt Komplexe Geometrie. }

\thanks{
\mbox{\hspace{11pt}}{\it 1991 Mathematics Subject Classification}:
14R05, 14R20, 14J50.\\
\mbox{\hspace{11pt}}{\it Key words}: $\C^*$-action, $\C_+$-action,
affine surface}

\dedicatory{Dedicated to Masayoshi Miyanishi}
\begin{document}

\begin{abstract}
Following an approach of Dolgachev, Pinkham and Demazure, we
classified in \cite{FlZa1} normal affine surfaces with hyperbolic
$\C^{*}$-actions in terms of pairs of $\Q$-divisors $(D_+,D_-)$ on
a smooth affine curve. In the present paper we show how to obtain
from this description a natural  equivariant completion of these
$\C^*$-surfaces. Using elementary transformations we deduce also
natural completions for which the boundary divisor is a
standard graph in the sense of \cite{FKZ} and show in certain
cases their uniqueness. This description is especially precise in
the case of normal affine surfaces completable by a zigzag i.e.,
by a linear chain of smooth rational curves. As an application we
classify all zigzags that appear as boundaries of smooth or normal
$\C^*$-surfaces.
\end{abstract}

\maketitle

\tableofcontents

\section{Introduction}

An irreducible normal affine surface $X=\Spec A$ endowed with an
effective $\C^*$-action will be called a {\it $\C^*$-surface}.
In the {\it elliptic case} the action possesses an
attractive or repulsive fixed point and in the {\it parabolic case}
an attractive or repulsive curve consisting of fixed points. A
simple and convenient description for these surfaces, based on
the fact that the
$\C^*$-action corresponds to a grading of the coordinate ring $A$
of $X$, was
elaborated by Dolgachev, Pinkham and Demazure, so it was called in
\cite[I]{FlZa1} a {\it DPD-presentation}. Namely, in the elliptic
case our surface is represented as $$X=\Spec A\qquad
\mbox{with}\quad A=\bigoplus_{k\ge 0} H^0(C,\cO_C(\lfloor kD
\rfloor))\cdot u^k\,,$$
where $u$ is an indeterminate, $D$ is an ample
$\Q$-divisor on a smooth projective curve $C$ and $\lfloor
kD\rfloor$ denotes the integral part. The curve $C=\Proj A$ is
then the orbit space of the $\C^*$-action on the complement of its
unique fixed point in $X$. Likewise, in the parabolic case
$$
X=\Spec A_0[D]\qquad \mbox{with}\quad
A_0[D] = \bigoplus_{k\ge 0} H^0(C,\cO(\lfloor kD \rfloor))\cdot
u^k\,,$$
where now $D$ is a $\Q$-divisor on a smooth affine curve
$C=\Spec A_0$, which again is the orbit space of our $\C^*$-action
on the complement of its fixed point set in $X$.

All other $\C^*$-surfaces $X$ are {\it hyperbolic}.
Their fixed points are all isolated, attractive in one and
repulsive in the other direction. Any
such surface is isomorphic to
$$
\Spec A_0[D_+,D_-]\qquad\mbox{with}\quad
A_0[D_+,D_-]:=A_0[D_+]\oplus_{A_0} A_0[D_-]\,
$$
where $D_\pm$ is a pair of $\Q$-divisors on a normal affine curve
$C=\Spec A_0$ with $D_++D_-\le 0$ \cite[I]{FlZa1}.

In this paper we are mainly interested in an explicit description
of the completions of such $\C^*$-surfaces.
One of the main results is contained in section 3,
where we describe a canonical equivariant completion of a
hyperbolic $\C^*$-surface in terms of the divisors $D_\pm$,
see for instance Corollary \ref{boundarydivisor} for
the dual graph
of its boundary divisor.
We also treat in brief the case of elliptic and parabolic surfaces,
see Section 3.4.

In \cite{FKZ}, Corollary 3.36 we have shown that any normal affine
surface $V$ admits a completion for which the dual graph of the
boundary is standard (see \ref{standardlist}). Given a DPD
presentation of a $\C^*$-surface $V$, the results of Section 3
provide an explicit equivariant standard completion $\bar V_{\rm
st}$ of $V$. More generally, in Section 2 we investigate the
question as to when such equivariant standard completions can be
found for actions of an arbitrary algebraic group $G$. We show
that this is indeed possible for normal affine $G$-surfaces $V$
except for
$$
\PP^2\backslash Q, \quad \PP^1\times\PP^1\backslash \Delta,
\quad V_{d,1}\,,
$$
where $Q$ is a non-singular quadric in $\PP^2$, $\Delta$ is the
diagonal in $\PP^1\times\PP^1$  and $V_{d,1}$, $d\ge 1$, are
the Veronese surfaces, see Theorem \ref{equivariant.5}.
Moreover, equivariant standard completions always exist if $G$ is
a torus. We also deduce their uniqueness in certain cases, see
Theorem \ref{unimain}.

In this paper we study mostly $\C^*$-actions on Gizatullin
surfaces.  By a {\it Gizatullin surface} we mean a normal affine
surface completable by a {\em zigzag} that is, a simple normal
crossing divisor $D$ with rational components and a linear dual
graph $\Gamma_D$. These surfaces are remarkable by a variety of
reasons. By a theorem of Gizatullin \cite[Theorems 2 and 3]{Gi}
(see also \cite{Be, BML}, and \cite{Du1} for the non-smooth case),
the automorphism group $\Aut (X)$ of a normal affine surface $X$
has an open orbit with a finite complement in $X$ if and only if
either $X\cong \C^*\times \C^*$ or $X$ is a Gizatullin
surface. The automorphism groups of Gizatullin surfaces were
further studied in \cite{DaGi}. Like in the case of $X=\A^2_\C$,
such a group has a natural structure of an amalgamated free
product.

These surfaces can also be characterized by the Makar-Limanov
invariant: a normal affine surface $X=\Spec A$ different from
$\A_\C^1\times\C^*$ is Gizatullin if and only if its Makar-Limanov
invariant is trivial that is, $\ML(X):=\bigcap \ker \p=\C\,,$
where the intersection is taken over all locally nilpotent
derivations of $A$. Among the hyperbolic $\C^*$-surfaces $X=\Spec
A_0[D_+,D_-]$ the Gizatullin ones are characterized by the
property that each of the fractional parts $\{D_\pm\}=D_\pm -
\lfloor D_\pm\rfloor$ is either zero or supported at one point
$\{p_\pm\}$, see \cite[II]{FlZa1}.

In Theorem \ref{thmzigzag}(a) we show that an
arbitrary ample zigzag  can be realized  as a boundary divisor of
a Gizatullin $\C^*$-surface and even a toric one.
However, not every such zigzag appears as the boundary divisor of
a {\em smooth} $\C^*$-surface. More precisely we give in
\ref{thmzigzag}-\ref{corzigzag} a numerical criterion as to when a
zigzag can be the boundary divisor of  a smooth Gizatullin
$\C^*$-surface. Using this criterion we can exhibit many smooth
Gizatullin  surfaces which do not admit any $\C^*$-action, see
Corollary
\ref{noaction}. We note that every $\Q$-acyclic Gizatullin
surface\footnote{That is $H_i(X,\Q)=0\quad\forall i>0$.} is a
$\C^*$-surface \cite[II.5.10]{Du2}. The latter class was studied
e.g., in \cite{DaiRu, MaMi1, Du2}.

Finally, in \ref{complet.15} we investigate $\C^*$-actions on
Danilov-Gizatullin surfaces, by which we mean complements
$\Sigma_n\setminus S$ of an ample section $S$ in a Hirzebruch
surface $\Sigma_n$. By a theorem of Danilov-Gizatullin
\cite{DaGi}, the isomorphism class of such a surface $V_{k+1}$
depends only on the self-intersection number $S^2=k+1> n$. In
particular it does not depend on $n$ and is stable under
deformations of $S$ inside $\Sigma_n$. According to Peter
Russell\footnote{An oral communication. We are grateful to Peter
Russell for generously sharing results from unpublished notes
\cite{CNR}.}, given any natural $k$ there are exactly $k$ pairwise
non-conjugated $\C^*$-actions on $V_{k+1}$. We give another proof
of this result using our DPD-presentations. In a forthcoming paper
we will show that a Gizatullin surface which possesses at least 2
non-conjugated $\C^*$-actions is isomorphic to a
Danilov-Gizatullin surface.

\section{Equivariant completions of affine $G$-surfaces}

\subsection{Equivariant completions}
\bsit\label{equivariant.1} By the Kambayashi-Mumford-Sumihiro
theorem (see \cite{Su}), any algebraic variety $X$ equipped with
an action of a connected algebraic group $G$ admits an equivariant
completion. For normal affine varieties this is true even without
the connectedness assumption. Indeed, if $X=\Spec A$ is an affine
$G$-variety then any $\C$-linear subspace of finite dimension of
$A$ is contained in a $G$-invariant one. Choosing an initial
$\C$-linear subspace which contains a set of algebra generators of
$A$ yields a $G$-invariant finite dimensional subspace
$E\subseteq
A$ such that the induced map gives an equivariant embedding
$X\hookrightarrow \A_\C^N$. Letting
$$
\A_\C^N\stackrel{\simeq}\longrightarrow
\A_\C^N\times\{1\}\subseteq \A_\C^N\times \A_\C^1$$
be a natural
embedding, where $G$ act on the second factor trivially, we get a
$G$-action on $\PP^N$. The closure $\bar X$ of $X$ in $\PP^N$ is
then an equivariant completion.

If $\dim\,X=2$ then an equivariant resolution of singularities of
such a completion can be obtained as follows. By a theorem of
Zariski \cite{Zar}, a resolution of singularities of $\bar X$ can
be achieved via a sequence of normalizations and blowups of points
i.e., of maximal ideals. Since both these operations are
equivariant, this yields an equivariant resolution. Moreover, the
minimal resolution dominated by this equivariant one is
equivariant too, provided that $G$ is connected and so  stabilizes
every component of the exceptional divisor. This is based on the
following  well known lemma, see e.g., Lemma 7 in \cite[I, \S
7]{DaGi}. \esit

\blem\label{equivariant.2} Let $X$ be a normal algebraic surface
with an action of an algebraic group $G$.
\begin{enumerate}[(a)]
\item Given a
contractible $G$-invariant complete curve $C$ in $X$, the action
of $G$ descends to the contraction $X/C$.
\item The action of $G$ lifts
to the blowup of $X$ in any fixed point of $G$.
\end{enumerate}
\elem

In the following, by an {\em NC completion} of a normal
algebraic surface $V$ we mean a pair $(X,D)$ such that $X$ is a
normal complete algebraic surface, $D$ is a normal crossing
divisor contained in the regular part $X_\reg$ and $V=X\backslash
D$. We call this an {\em SNC completion} if moreover $D$ has only
simple normal crossings.

The  considerations above lead to the following well known result.

\bprop\label{pSNC}
(a) A normal affine algebraic surface $V$ with an action of an
algebraic
group admits an equivariant SNC completion $(X,D)$.

(b) An arbitrary normal algebraic surface $V$ with an action of a
connected algebraic group admits an equivariant SNC completion
$(X,D)$.

(c)
Any two equivariant SNC completions $(X_i, D_i)$ of
$V$, $i=1,2$, are equivariantly dominated by a third one $(X, D)$.
\eprop

\bsit Let $\G$ be a weighted graph. We recall (see Definitions 2.3
and 2.8 in \cite{FKZ}) that an {\em inner} blowup $\G'\to\G$ is
one performed in an edge of $\Gamma$, and that an  {\em
admissible} blowup is one that is inner or performed in an end
vertex of $\Gamma$. Moreover a blowdown $\G\to \G''$ is said to be
admissible if its inverse is so. A  birational transformation of
graphs is a sequence of blowups and blowdowns. Given such a
sequence \be\la{combi}\bdi \gamma:\quad
\G=\G_0&\rDotsto^{\gamma_1}& \G_1 & \rDotsto^{\gamma_2} &\cdots &
\rDotsto^{\gamma_n}& \G_n=\G'\edi\,, \ee we call it {\em
admissible} if every $\gamma_i$ is so, and {\em inner} if every
step is an admissible  blowdown or an inner blowup.
 \esit

\bdefi\label{reconstruction} Given two NC completions
$(X,D)$, $(X',D')$ of a normal algebraic surface $V=X\backslash
D=X'\backslash D'$, by a {\em birational map} $\psi:
(X,D)\dashto(X',D')$ we mean a birational map  $X\dashto X'$
inducing the identity on $V$. Such a map can be decomposed
into a sequence of blowups and blowdowns \be\la{combi2}\bdi
\tgamma:\quad (X,D)=(X_0,D_0)&\rDotsto^{\tgamma_1}& (X_1,D_1) &
\rDotsto^{\tgamma_2} &\cdots & \rDotsto^{\tgamma_n}&
(X_n,D_n)=(X',D')\edi\,, \ee where (i) $X_{i+1}$ is a blowdown or
a blowup of $X_i$ taking place in the total transform $D_i$ of $D$
in $X_i$ and (ii) $D'$ is the total transform of $D$. Clearly
$\tgamma$ will induce a birational map $\gamma$ as in
(\ref{combi}) of the dual graphs $\G_i$ of $D_i$.

A birational map $\psi: (X,D)\to(X',D')$ will be called {\em
inner} or {\em admissible} if $\gamma$ has the respective property
for a suitable factorization $\tgamma$ as above. If $X$ is
equipped with an action of an algebraic group $G$ leaving $D$
invariant, then we call $\psi$ or the sequence $\tgamma$ {\em
$G$-equivariant} if they are compatible with the action of $G$.
\edefi

The following observation will be useful.

\bprop\label{equivariant.4} Let $G$ be a connected algebraic group
acting on a normal algebraic surface $V$ and let $(X,D)$ be an
equivariant NC completion of $V$. Assume that $\gamma:\G\dashto
\G'$ is a birational transformation of the dual graph $\Gamma$ of
$D$ as in (\ref{combi}) that blows down at most vertices of $\G$
corresponding to rational components of $D$. Then there is a
sequence of equivariant birational maps $\tgamma:(X,D)\dashto
(X',D')$ as in (\ref{combi2}) inducing $\gamma$ on the dual graphs
of $D$, $D'$ in each of the following cases.
\begin{enumerate}[(i)]
\item $\gamma$ is
inner.
\item
$G=\T=(\C^*)^n$ is a torus and $\gamma$ is admissible.
\end{enumerate}
\eprop

\bproof
(i) is immediate from Lemma \ref{equivariant.2}, and (ii)
follows as well since an action of a torus on the projective
line has at least 2 fix points.
\eproof From this Proposition we can
deduce the following corollaries.

\bcor\label{equivariant.05}
For a normal surface $V$ with an action of a
connected algebraic group $G$ the following hold.

(a)  $V$ admits a minimal equivariant NC completion $(X,D)$, i.e.\
$D$ contains no at most linear\footnote{i.e. such that the
degree of the corresponding vertex in the dual graph of $D$ is
$\le 2$.} rational $(-1)$-curve.

(b) If moreover $G=\T$ is a torus and $(X,D)$ and $(X',D')$ are
two minimal equivariant NC completions of $V$ then there is an
equivariant admissible birational map $\psi:(X,D)\dashto (X',D')$.
\ecor

\bproof (a) is an immediate consequence of Proposition
\ref{equivariant.4}. If all irreducible components of $D$ (and
then also of $D'$) are rational curves then (b) follows from
Propositions 2.9 in \cite{FKZ} and \ref{equivariant.4}. In the
general case we proceed as follows. If  $v$ is a vertex of the
dual graph $\G$  of $D$ corresponding to a non-rational curve then
we add a simple loop at $v$. This procedure results in a new
minimal graph $\tilde \G$ in which the vertices corresponding to
non-rational curves become branching points. In the same way we
obtain from the dual graph $\G'$ of $D'$ a graph $\tilde\G'$ that
is birationally equivalent to $\tilde \G$. According to
Proposition 2.9 in \cite{FKZ} $\tilde\G'$ can be obtained from
$\tilde\G$ by an admissible birational transformation. Omitting at
each step the simple loops just added results in  an admissible
birational transformation of $\G$ into $\G'$.
Applying Proposition
\ref{equivariant.4} the assertion follows. \eproof

\subsection{Standard and semistandard completions}
We use below
the notions of standard and semistandard graphs as introduced in
\cite[Definition 2.13]{FKZ}. For the convenience of the reader
we recall some of the notations from \cite{FKZ}.

\bsit\label{standardlist} Since the dual weighted graph of a
divisor on an algebraic surface satisfies the Hodge index theorem
we restrict in the sequel to graphs whose intersection form
has at most one positive eigenvalue. Following the notations in
\cite{FKZ} we use the abbreviation \be\label{notation}
[[w_1,\ldots, w_n]]\quad:=\qquad \co{w_1}\lin \co{w_2}
\lin\ldots\lin\co{w_n}\quad, \ee and $((w_1,\ldots, w_n))$ will
denote the circular standard graphs obtained from this by
connecting the first and last vertex by an additional edge.

A graph  $[[w_1,\ldots, w_n]]$ (or $((w_1,\ldots, w_n))$\,)
will be called a (circular) {\em zigzag} if its intersection
form has at most one positive eigenvalue. According to
\cite[Lemma 2.17 and
Proposition 4.13]{FKZ} the standard zigzags are
\be
\label{standard} [[0]]\,, \quad [[0,0,0]]
\,\quad\mbox{and}\quad
  [[0,0,w_1,\ldots,w_n]],\quad \mbox{where}
\quad n\ge 0,\, w_j\le-2\,\,\,\forall j, \ee
and the circular standard zigzags
\be
\label{cstandard}
((0_a,w)),\quad ((0_b,-1,-1))\quad \mbox{and}\quad
((0_{b}, w_1,\ldots, w_n))\,,
\ee
where $0\le a\le 3$, $w\le 0$, $b\in\{0,2\}$ and
$w_i\le -2$ $\forall i$.
In geometry there also appear naturally
semistandard zigzags, where we have additionally the possibilities
\be\label{mstandard}  [[0,w_1,\ldots,w_n]]\,, \,\,\,
[[0,w_1,0]], \qquad\mbox{where }
n\ge 0\quad\mbox{and}\quad w_j\le-2\,\,\,\forall j\,,
\ee
see \cite[Lemma 2.17]{FKZ}.

We notice that a standard zigzag $[[0,0,w_1,\ldots,w_n]]$ is
unique in its birational class up to reversion
\be\la{reverse}
[[0,0,w_1,\ldots,w_n]]\rightsquigarrow [[0,0,w_n,\ldots,w_1]]\,,
\ee
and the circular standard zigzag $((0_{b}, w_1,\ldots, w_n))$ is
unique up to reversion and a cyclic permutation
$$((0_{b}, w_1,\ldots, w_n))\rightsquigarrow ((0_{b},
w_{q-1},\ldots,w_n,\,w_1,\ldots, w_q))\,.$$  The other standard
zigzags are unique, see Corollary 3.33 in \cite{FKZ}.

In the following an NC divisor $D$ with dual graph $\G$ on
an algebraic surface will be called {\it standard} or {\it
semistandard} if all connected components of $\G\ominus
(B(\G)\cup S)$ have this property, where $B(\G)$ is the set of all
branching points of $\G$ and $S$ is the set of vertices
corresponding to non-rational curves. Similarly, a completion
$(X,D)$ of an open algebraic surface is said to be (semi-)standard
if $D$ is so. \esit

The next result  is an analogue of Theorem 7 in \cite[I]{DaGi}
which says that any algebraic group action on an affine surface
admitting a standard  completion (in the sense of \cite{DaGi}),
admits also an equivariant standard completion. However note that
our standard zigzags form a narrow subclass of those in
\cite[I]{DaGi}.

\bthm\label{equivariant.5} (a) Every normal affine surface $V$ with
an action of a connected algebraic group $G$ admits an
equivariant semistandard NC completion $(X,D)$ unless $X$ is one
of the surfaces
$$
\PP^2\backslash Q, \quad \PP^1\times\PP^1\backslash \Delta,
\quad V_{d,1}\,,
$$
where $Q$ is a non-singular quadric in $\PP^2$, $\Delta$ is the
diagonal in $\PP^1\times\PP^1$ and $V_{d,1}$, $d\ge 1$, are
the Veronese surfaces\footnote{See e.g. Lemma \ref{torics}(a)
below.}.

(b) If $G=\T$ is a torus and $V$ is an arbitrary normal surface
then there is an equivariant standard completion  $(X,D)$. \ethm

\bproof

Let $(Y,E)$ be an equivariant NC completion of $V$. Let us
first suppose that $E$ is not an irreducible smooth rational curve
so that the dual graph $\G$ of $E$ is not reduced to a point. If
all components of $E$ are rational then by Theorem 2.15 in
\cite{FKZ} $\G$  can be transformed into a semistandard graph by
an inner birational transformation and even into a standard one
by an
admissible transformation.  Thus  both claims follow now from
Proposition
\ref{equivariant.4}. If some of the components are not rational,
then  as in the proof of Corollary \ref{equivariant.05} we
can add
to $\G$ simple loops so that the vertices corresponding to
non-rational curves become branching points. Arguing as before
the
result also follows in this case.

Assume further that $E$ is a smooth irreducible rational curve. If
the group $G$ is solvable then there is a fixed point of $G$ on
$E$, and blowing it up successively we can transform $E$ into a
chain $[[0, -1,-2,\ldots, -2]]$, see \cite[Remark 2.14(1)]{FKZ}.
Since this chain can be transformed into a semistandard (standard)
one by an equivariant inner (admissible) elementary transformation
the result follows also in this case.

Finally, if $G$ is not solvable then it contains  a subgroup
isomorphic to $\SL_2(\C)$ or $\PGL_2(\C)$. Using the theorem
of Gizatullin and Popov (see Proposition 4.14 in \cite{FlZa2}
and the references therein) our surface is one of the list
above. \eproof

\brems 1. As the proof shows, (a) holds for an arbitrary normal
algebraic surface $V$ provided that $G$ is solvable or $V$ admits
an equivariant NC completion $(Y,E)$ such that the dual graph of
$E$ is not reduced to a point.

2. We cannot expect in general to obtain an equivariant standard
completion for a solvable group, because there could be not enough
fixed points to perform outer equivariant elementary
transformations as required to get a standard form. For instance,
the group $G$ of all projective transformations of $\PP^2$ which
stabilize a line $D$ and a point $A\in D$ is solvable and has the
only fixed point $A$. There exists an equivariant completion of
$\A^2_\C=\PP^2\setminus D$ with semistandard dual graph
$[[0,-2]]$, but it is impossible to get such a completion with
standard dual graph $[[0,0]]$. \erems

Next we address the question of uniqueness of (semi-)standard
completions. We recall shortly the notion of elementary
transformations. Given a linear 0-vertex $v$ of $\G$, so that $\G$
contains $L=[[w,0,w']]$ we consider the birational map of $\G$
given by
\begin{equation}\label{etr1}
[[w-1,0, w'+1]]\dashrightarrow [[w-1, -1,-1,w']] \lto [[w,0,w']].
\end{equation}
on $L$, which is the identity on $\G\ominus L$. Similarly, if
$v\in\G$ is an end vertex so that $\G$ contains $L'=[[w,0]]$,
we consider the birational map of $\G$ given on $L'$ by
\begin{equation}\label{etr2}
[[w-1,0]]\dashrightarrow [[w-1, -1,-1]] \lto [[w,0]]\,.
\end{equation}
These transformations as well as their inverses are called
{\em elementary transformations} of $\G$.

Similarly, given a completion $(X,D)$ of a normal surface $V$ we
can define elementary transformations at any point of a component
$C_i\cong \PP^1$ of $D$ of selfintersection 0 that corresponds to
an at most linear vertex of the dual graph of $D$.

\bprop Let $G$ be a connected algebraic group acting on a normal
algebraic surface $V$. If $(X_1,D_1)$ and $(X_2,D_2)$ are
equivariant semistandard NC completions of $V$, then $(X_2,D_2)$
can be obtained from $(X_1,D_1)$ by a sequence of equivariant
elementary transformations of the boundary. \eprop

\bproof Let us first assume that the irreducible components of
$D_1$ and $D_2$ are all rational. By Proposition \ref{pSNC} there
is an equivariant NC completion $(X,D)$ of $V$ dominating
$(X_i,D_i)$ for $i=1,2$. If $\G$, $\G_1$ and $\G_2$ are the
respective dual graphs of $D$, $D_1$ and $D_2$ then $\G$ dominates
$\G_1$ and $\Gamma_2$. By Theorem 3.1 in \cite{FKZ} we can
transform $\G_1$ into $\G_2$ by a sequence of elementary
transformations such that every step is dominated by some inner
blowup of $\G$. Using Proposition 3.34 from \cite{FKZ} this
gives a unique sequence of elementary transformations transforming
$(X_1,D_1)$ into $(X_2,D_2)$ such that every step is dominated by
an inner blowup , say $(X',D')$, of $(X,D)$. Since by Lemma
\ref{equivariant.2} the action of $G$ lifts naturally to $(X',
D')$ and $G$ also acts on any blowdown of the boundary $D'$, the
result follows in this case.

In the general case we can again add simple loops at the vertices
of $\G_1$, $\G_2$ and $\G$ as in the proof Corollary
\ref{equivariant.05}. Arguing as before the result follows
also in this case.
\eproof

\subsection{Uniqueness of standard completions}
In general, standard equivariant completions even of
$\C^*$-surfaces are by no means unique. Let us give two examples.

\bexa \la{exuni} 1. Given a Gizatullin $\C^*$-surface $V$ and an
equivariant standard completion $(V,D)$ we can reverse the
boundary zigzag $D$ as in (\ref{reverse}) by a sequence of inner
elementary transformations. This leads to another equivariant
standard completion, which usually is not isomorphic to
the given one.

2. The affine plane
$\A^2$ endowed with the $\C^*$-action $t.(x,y)=(tx,ty)$
can be equivariantly completed by $\PP^1\times\PP^1$.
The dual graph of the
boundary divisor is the standard zigzag
$[[0,0]]$ consisting of the curves, say $C_0$ and $C_1$.
Blowing up the intersection point $C_0\cap C_1$ and blowing down
$C_1$ gives a component, say $E$ that is pointwise fixed by $\C^*$.  Performing an outer
blowup of $E$ in a point different from the contraction of $C_1$,
and  then blowing down $E$, we arrive at a new equivariant
completion of $\A^2$ by a standard zigzag as before. However, the
equivariant completions of $\A^2$ obtained in this way are not
equivariantly isomorphic, although both of them are isomorphic to
$\PP^1\times\PP^1$ and the boundary zigzags are the same.
\eexa

The main result of this section is the following
uniqueness theorem.

\bthm\label{unimain} \begin{enumerate}[(a)] \item A non-toric
Gizatullin $\C^*$-surface $V$ has a unique standard completion up
to reversing the boundary zigzag. More precisely, any two such
completions $(\bV_\st, D_\st)$ and $(\bV'_\st,D'_\st)$ are
isomorphic or obtained from each other by reversing the boundary
zigzag.

\item A normal affine toric surface $V$ has a unique standard
completion up to reversing the boundary zigzag unless $V$ is one
of the surfaces $\A\times\C^*$ or $\C^*\times\C^*$.
\end{enumerate}
\ethm

The assertions (a), (b) of the theorem will be shown in
\ref{proofa} and \ref{proofb} below, respectively. We need a few
preparations.

\bdefi\label{6.8} Let $V$ be a normal surface with an action of an
algebraic group $G$. A curve of fixed points of $G$ in $V$ will be
called {\em $G$-parabolic}, or simply {\em parabolic} if $G$ is
clear from the context.
\edefi

The following lemma is well known. For the sake of completeness we
provide a simple argument.

\blem\label{before01} Let the $2$-torus $\T$ act on $V_0\cong
\C^*\times\C^*$ with an open orbit, and let $(\bV,D_0)$ be an
equivariant smooth completion of $V_0$ by an SNC divisor $D_0$.
Then $D_0$ is a cycle of rational curves without
$\T$-parabolic components. \elem

\bproof As follows e.g., from Luna's \'Etale Slice Theorem, for
any regular action of an algebraic reductive group with an open
orbit, the fixed point set is finite. (In the toric case there is
an easy direct argument; cf. \cite{Su}.) Hence $\bV$ cannot
contain $\T$-parabolic curves.

The surface $V_0\simeq \C^*\times\C^*$ admits an equivariant
completion  $(\PP^1\times\PP^1, Z_0)$ by a cycle $Z_0$ consisting
of $4$ rational curves. Thus there is an equivariant birational
transformation $\gamma: D_0\dashrightarrow Z_0$. We claim that
$\gamma$ is inner, so at each step of this transformation the
boundary divisor remains a cycle of rational curves, as required.
Indeed, this follows by induction on the length of $\gamma$, using
the fact that $\gamma$ can blow up only isolated fixed points of
$\T$ on the boundary, which are double points of the boundary
cycle by the inductive hypothesis.
\eproof

\bsit\la{proofb}{\em Proof of Theorem \ref{unimain}(b)}.
If $V$ is not one of the surfaces  $\A\times\C^*$ or
$\C^*\times\C^*$,
then the boundary zigzag of a standard completion
$(\bV_\st, D_\st)$ is not equal to $[[0,0,0]]$ and
is not circular.
Comparing with the list in
(\ref{standard}) the dual graph of $D_\st$ is of the form
$[[0,0,w_2,\ldots, w_n]]$ with $w_i\le -2$ for all $i$
and $n\ge 1$.
Given another standard completion $(\bV'_\st, D_\st')$
there is an  equivariant domination $(Y,E)$ of these completions.
By Lemma \ref{before01} $E$ is again a zigzag and so by
Proposition 3.4
in \cite{FKZ} $\bV_\st=\bV'_\st$, or $D_\st'$ is obtained
by reversing
the zigzag $D_\st$.
\esit

We now embark on the proof of the more difficult part (a) of
Theorem \ref{unimain}. Let us first fix some notations.

\bsit\la{unisit} Let $V$ be a non-toric Gizatullin surface and
$(\bV_\st, D_\st)$ a completion of $V$ by a standard zigzag
$[[0,0,w_2,\ldots, w_n]]$ with $w_i\le -2$ $\forall i$ and $n\ge
2$. We let
$$
D_\st=C_0+\ldots+C_n\,,
$$
where the components are numbered according to the weights in
the sequence \newline $[[0,0,w_2,\ldots, w_n]]$. We also consider
the minimal resolutions of singularities $V'$, $(\tV_\st,D_\st)$
of $V$  and $(\bV_\st, D_\st)$, respectively.

Since $C_0^2=C_1^2=0$ the
linear systems $|C_0|$ and $|C_1|$ define a morphism
$\Phi=\Phi_0\times\Phi_1:\tV_{\rm st}\to \PP^1\times\PP^1$ with
$\Phi_i=\Phi_{|C_i|}$, $i=0,1$. We notice that $C_1$ is a section
of $\Phi_0$ and so the restriction $\Phi_0|V':V'\to\PP^1$ is an
$\A^1$-fibration. We can choose the coordinates in such a way that
$$C_0=\Phi_0^{-1}(\infty)\,, \qquad \Phi(C_1)= \PP^1\times
\{\infty\} \quad\mbox{and}\quad C_2\cup\ldots \cup C_n\subseteq
\Phi_0^{-1}(0) \,.$$ The divisor $D_\ext:=C_0\cup C_1\cup
\Phi_0^{-1}(0)$ is called the {\it extended divisor}. It
will be studied systematically in Section 5. \esit

\brem\label{unirem}
If $V$ carries a $\C^*$-action, then we can find equivariant
standard completions  $(\bV_\st, D_\st)$ and  $(\tV_\st, D_\st)$,
see Proposition \ref{equivariant.5}. Thus $\Phi$
will also be equivariant with a suitable
$\C^*$-action on $\PP^1\times\PP^1$.
\erem

\blem\label{unilem}
With the notation as in \ref{unisit},
$\Phi$ is birational and induces an isomorphism
$\tV_\st\backslash \Phi_0^{-1}(0)\cong
(\PP^1\backslash\{0\})\times \PP^1$.
In particular, $\Phi_0^{-1}(0)$ is the only
possible degenerate fiber of the
$\PP^1$-fibration $\Phi_0 : \tV_\st\to \PP^1$.
\elem

\bproof
Since by construction $\Phi^{-1}(\infty, \infty)=C_0\cap C_1$
consists of one point, the map is birational and so $\tV_\st$
is a blowup of $\PP^1\times\PP^1$. Because of $C_0^2=C_1^2=0$
no blowup can occur along $C_0\cup C_1$, whence $\Phi$ is an
isomorphism in a neighborhood of $C_0\cup C_1$.

Now assume that for some point $x\in (\PP^1\backslash\{0,
\infty\}) \times (\PP^1\backslash\{\infty\})$ the fibre
$\Phi^{-1}(x)$ is a curve. Then this curve meets neither
$C_0\cup C_1$ nor the divisor $D_\st\ominus C_0\ominus C_1$ since
by construction, the latter one  is contained in
$\Phi_0^{-1}(0)$. Thus $\Phi^{-1}(x)$ is contained in $V'$.
Since $V$ being affine does not contain complete curves
this is only possible if $\Phi^{-1}(x)$ is contained in the
exceptional divisor of $V'\to V$.  Because $\Phi^{-1}(x)$
contracts to a smooth point in $\PP^1\times\PP^1$ it must contain
a $(-1)$-curve,  which gives a contradiction since $V'$ is the
minimal resolution of $V$. \eproof

\blem\label{pextended.01} In the notation of \ref{unisit}, if for
some standard completion $(\bV_\st, D_\st)$ of a Gizatullin
surface $V$ the extended divisor $D_\ext$ is linear  then $V$ is
toric.

On the other hand, for any {\rm equivariant} standard completion
$(\bV_\st, D_\st)$ of a toric Gizatullin surface $V$ the extended
divisor $D_\ext$ is linear. \elem

\bproof Since $C_1^2=0$ both on $\tV_{\rm st}$ and on
$Q=\PP^1\times\PP^1$, no blowup is done under
$\Phi=(\Phi_0,\Phi_1): \tV_{\rm st}\to Q$ with center on $C_1$.
Thus we may assume that the center of the first blowup in $\Phi$
is the fixed point $(0,0)\in C_2\ominus C_1$ of the standard
$\T$-action on $Q$. We claim that this action lifts to $\tV_{\rm
st}$ stabilizing $V'$ and then descends to $V$.

Indeed, by Lemma \ref{unilem}, $\Phi_0^{-1}(0)$ is the only
possible degenerate fiber of $\Phi_0$. Thus, with $E$
the exceptional set of the minimal resolution $V'\to V$,
both $D_{\rm st},E$ are
disjoint subchains of the linear chain $D_{\rm ext}$.
Since $V$ is
affine, contracting $E$ every component of $D_{\rm ext}\ominus
(D_{\rm st}+ E)$ meets the image of $D_{\rm st}$. Since $D_{\rm
st}$ is connected it follows that there is exactly one such
component, say, $E_0$ which separates $D_{\rm st}$ and $E$.
Moreover since $D_{\rm st}$ and $E$ are both minimal, $E_0$
is the
only $(-1)$-curve in $D_{\rm ext}$.

Therefore all blowups in $\Phi|D_{\rm ext}: D_{\rm ext}\to
C_0+C_1+C_2$ are inner except for the first one. Hence the
standard torus action on $Q$ lifts through $\Phi$ to $\tV_{\rm
st}$ leaving $D_{\rm ext}$ stable.  It stabilizes as well
$D_{\rm st},E$ and $V'=\tV_{\rm st}\setminus D_{\rm st}$ and so by
Lemma \ref{equivariant.2} descends to $V$  with an open orbit.
Thus indeed $V$ is toric.

As for the converse, note that by Lemma \ref{before01} $D_{\rm
ext}$ is part of a cycle of rational curves. Hence being connected
and simply connected, it is a linear chain.\eproof

\blem\label{parabolic} Assume that $(\bV_\st, D_\st)$ is an
equivariant completion of a normal affine $\C^*$-surface $V$.
With
the notation as in \ref{unisit}, if $V$ is non-toric then one of
the curves $C_2,\ldots, C_n$ is parabolic. \elem

\bproof We note that the fiber $\Phi_0^{-1}(0)$ is invariant
under
the $\C^*$-action. Since $V$ is non-toric, by Proposition
\ref{pextended.01} the dual graph of $D_{\rm ext}=C_0\cup
C_1\cup\Phi_0^{-1}(0)$ contains a branching point $C_k$, $k\ge 2$.
Thus the $\C^*$-action has at least 3 fixed points on this
component $C_k$ which is then parabolic, as needed.\eproof

\bsit\la{proofa}{\em Proof of Theorem \ref{unimain}(a)}. Let
$(\bV_\st, D_\st)$ be an equivariant standard completion of $V$.
With the notations as in \ref{unisit}, by Lemma \ref{parabolic}
there is a parabolic component, say, $C_{s+1}$ in $D_\st$. After
moving  the 2 zero weights in the zigzag via a sequence of inner
elementary transformations to the components $C_{s}$ and $C_{s+1}$
we get a new equivariant completion $(\tV,D)$ of $V$. Note that
moving these zeros the curve $C_{s+1}$ is not blown down, and that
the inverse transformation $D\dashto D_\st$ is as well inner, cf.\
Lemma 2.12 in \cite{FKZ}. The linear system $|C_s|$ gives a
morphism $\psi:\tV\to \PP^1$ equivariant with respect to a
suitable $\C^*$-action on $\PP^1$, where $\psi(C_s) =\infty$.  The
curves $C_{s\pm 1}$ being disjoint sections of $\psi$ and
$C_{s+1}$ being parabolic, $\psi$ is the orbit map. We let $\bV$
be the surface obtained from $\tV$ by contracting all curves in
$D$ besides $C_{s\pm 1}$ and $C_s$. Obviously $\tV$ is then the
minimal resolution of the singularities of $\bV$ sitting on the
boundary.

Given a second equivariant standard completion $(V'_\st, D'_\st)$,
with the same procedure we get surfaces $\tV'$ and $\bV'$
fibered equivariantly over $\PP^1$. As before $\bV'$ is a
completion of $V$ by three curves $C'_{t-1}$, $C'_t$ and
$C_{t+1}'$ so that $C'_{t-1}$ and $C'_{t+1}$ are sections
of the $\PP^1$-fibration and $C'_t$ is the fiber over $\infty$.
The identity map on $V$ extends to an equivariant
birational map $h:\bV \dashto \bV'$ compatible with
the orbit maps to $\PP^1$. In particular, $h$ respects sections
of the $\PP^1$-fibrations and so $C_{s+1}$ is the proper
transforms of one of the sections $C'_{t+ 1}$ or $C'_{t-1}$
in $\bV'$, and similarly for $C_{s-1}$.
Performing, if necessary, elementary transformations at the fiber
$C_s$ we may also assume that $C_s$ is the proper transform of
$C'_t$.

Now $h$ defines a biregular map on the complements of discrete
sets, so by Zariski's main theorem, it is everywhere regular and
an isomorphism. This isomorphism  lifts to the minimal resolutions
of singularities giving an equivariant isomorphism $\tilde h:
(\tV,D)\to (\tV',D')$. Since $(\tV,D)\dashto (\bV_\st, D_\st)$
and
$(\bV'_\st, D'_\st)\dashto (\tV',D')\cong (\tV, D)$ are both
composed of inner elementary transformations it follows that
$D'_\st\dashto D_\st$ is as well  inner. Thus using
Proposition 3.4 in \cite{FKZ}, either $D_\st=D'_\st$, or
$D_\st $ is the reversion of $D'_\st$, proving (b). \qed\esit

\brem For an arbitrary normal affine $\C^*$-surface $V$ the dual
graph of a standard equivariant completion can be easily deduced
from the description in Corollary \ref{boundarydivisor}(b). It is
easy to see that, if the surface is not a Gizatullin one, it
admits in general  many different equivariant standard
completions. \erem

\section{Equivariant completions of $\C^*$-surfaces}

\subsection{Generalities}
\bsit\label{complet.1} For an arbitrary normal compact complex
surface $X$, there is a $\Q$-valued intersection theory for
divisors on $X$ (see \cite[\S II.4]{Mu}, \cite{Sa}). This is a
pairing
$$
\Div_\Q(X)\times \Div_\Q(X)\to\Q,\quad (D_1,D_2)\mapsto
D_1.D_2\in \Q\,,
$$
sharing the usual properties of intersections on smooth surfaces:

\begin{enumerate}[1.]
\item The pairing is bilinear.
\item The projection formula with respect to proper mappings
$f:X\to Y$
of normal surfaces holds:
$$
f^*(D).E=D.f_*(E) \quad\mbox{for }D\in \Div_\Q(Y) \mbox{ and
}E\in\Div_\Q(X)\,.
$$
\item The adjunction formula holds, i.e.\ if $C\subseteq X$ is an
integral curve
and $D$ is a Cartier divisor on $X$
then $C.D=\deg_C(\cO_X(D)|C)$.
\end{enumerate}
\esit

For a sequence of real numbers $k_0,\ldots, k_n$ with
$k_0,\ldots, k_{n-1}\ge 2$ and $k_n\ge 1$ we let
$[k_0,\ldots,k_n]$ be the continued fraction defined inductively
via
$$
[k_0]=k_0\quad\mbox{and}\quad[k_0,\ldots,k_n]=
k_0-\frac{1}{[k_1,\ldots,k_n]}\,\,.
$$

\bprop\label{complet.2} Let $X$ be a normal surface and let
$C_0,C_1,\ldots, C_n$ be a chain of irreducible curves with
$C_{i-1}.C_i=1$ for $i=1,\ldots, n$ and $C_i.C_j=0$ for $i\ne j$
otherwise, and with dual graph
$$
\cou{C_0}{-k_0}\lin\cou{C_1}{-k_1}\ldots\quad \ldots
\cou{C_{n-1}}{-k_{n-1}}\lin\cou{C_n}{-k_n}\quad ,
$$
where $k_i=-C_i^2\ge 2$ $\forall i=1,\ldots,n$ (however, we allow
$X$ and the $C_i$ to be singular so that $k_i\in \Q$). Assume that
$C_1\cup\ldots\cup C_n$ can be contracted via a map $\pi:X\to X'$,
and let $C_0'=\pi(C_0)$ be the image of $C_0$ in $X'$. Then
\be\label{intfor} -C_0^{\prime
2}=[k_0,k_1,\ldots,k_n]=k_0-\frac{1}{[k_1,\ldots, k_n]}\,. \ee In
particular, in the case where $k_i\in\N$ $\forall i$ we have
$C_0^2=\lfloor C_0^{\prime 2}\rfloor$.
  \eprop

\bproof We write $\pi^*(C_0')=C_0+r_1C_1+\ldots+r_nC_n$. By the
projection formula
$$
\pi^*(C_0').C_0=C_0'^2 \quad\mbox{and}\quad
\pi^*(C_0').C_i=0\mbox{ for }1\le i\le n.
$$
This leads to the equalities
$$
C_0^{\prime 2}=C_0^2+r_1 \quad\mbox{and}\quad \frac{r_{i-1}}{r_i}
= k_i-\frac{r_{i+1}}{r_{i}}\mbox{ for }1\le i\le n\,,
$$
with the convention that $r_0=1$ and $r_{n+1}=0$. Hence by
induction
$$
\frac{r_{i-1}}{r_i}=[k_i,\ldots, k_n]\,.
$$
In particular, $\frac{r_0}{r_1}=[k_1,\ldots, k_n]$. As
$-C_0^{\prime 2}=-C_0^2-r_1$, (\ref{intfor}) follows. The
last assertion also follows as $C_0^2=C_0^{\prime 2}-r_1\in\Z$,
where $0<r_1=[k_1,\ldots, k_n]^{-1}<1$ by our assumption.
\eproof

\bexa\label{complet.3} Suppose that $X$ as in \ref{complet.2} is
smooth and that $C_1,\ldots, C_n$ is a chain of smooth
$(-2)$-curves. In this case $[2,\ldots, 2]=(n+1)/n$ and so,
$C_0^{\prime 2}=C_0^2+n/(n+1)$. For instance, if $C_0$ is a
$(-1)$-curve in $X$ then the self-intersection number of $C_0'$ is
$-1/(n+1)$. \eexa

\brem\label{complet.4} If the curves $C_1,\ldots, C_n$ as in
\ref{complet.2} above are smooth and sitting in the smooth locus
of $X$ then by a result of Grauert \cite{Gr} these curves can be
contracted in the category of normal analytic spaces, provided
that $k_i\ge 2$ $\forall i=1,\ldots,n$. However, in general $X'$
is not necessarily a scheme, see for instance \cite{Sch}. \erem

\blem\label{complet.5} Let $D\in\Div_\Q(C)$ be a $\Q$-divisor on
a smooth complete curve $C$ and let $\cO_C[D]$ be the sheaf of
$\cO_C$-algebras
$$
\cO_C[D]:=\bigoplus_{k\ge 0} \cO_C(\lfloor kD \rfloor)
\cdot u^k\,,$$
where $u$ is an indeterminate.
The (relative) spectrum $X=\Spec \cO_C[D]$ is then
a normal surface, and
the zero section $S\subseteq X$ corresponding to the projection
$\cO_C[D]\to \cO_C$ has selfintersection $-\deg (D)$. \elem

\bproof
Consider $d\in \N$ such that the divisor
$D'=dD$ is Cartier on $C\simeq S$. If $\zeta$ is a local
generator of
$\cO_C(dD))$ in a neighbourhood of a point $s \in S$ as an
$\cO_C$-module then $dS$ is given by the zeros of the local
section $\zeta u^d$ in $\cO_C[D]$.
Thus $dS$ is
Cartier on $X$. By adjunction $dS.S=\deg (dS|S)$. Under
the canonical identification $S\cong C$ we have $dS|S=-dD$, so
$dS^2=-d\deg (D)$, proving the lemma. \eproof

\subsection{Equivariant completions of hyperbolic $\C^*$-surfaces}

\bsit\label{complet.0} In this subsection $V$ denotes a hyperbolic
$\C^*$-surface. According to \cite[I]{FlZa1}, such a surface is
isomorphic to $\Spec A_0[D_+,D_-]$, where $D_\pm$ is a pair of
$\Q$-divisors on the normal affine curve $C=\Spec A_0$ with
$D_++D_-\le 0$. This means that
$$
A=A_{\le 0}\oplus_{A_0}A_{\ge 0}\subseteq K[u,u^{-1}]
$$
with $K=\Frac (A_0)$,
$$
A_{\ge 0}=\bigoplus_{i\ge 0} H^0(C,\cO_C(\lfloor
iD_+\rfloor))\cdot u^i \quad\mbox{and}\quad A_{\le
0}=\bigoplus_{i\le 0} H^0(C,\cO_C(\lfloor -iD_-\rfloor))\cdot
u^i\,.
$$
Our goal is to describe a canonical completion of such a
$\C^*$-surface $V$ in terms of the divisors $D_\pm$.\esit

\bsit\label{glue}
Let us consider the same pair of $\Q$-divisors
$D_\pm$ on the smooth completion $\bC$ of $C$.
Identifying
the function field $K=\Frac(\bC)$ with the constant sheaf $K$ on
$\bC$, we form the sheaf of $\cO_\bC$-algebras
$$
\cO_\bC[D_+,D_-]\subseteq K[u,u^{-1}]
$$
by defining it on affine open subsets as in \ref{complet.0}. The
resulting normal $\C^*$-surface $V_0=\Spec \cO_\bC[D_+,D_-]$
contains $V$ as an open subset and can be completed as follows.
\esit

\bprop\label{complet.6} There is a natural $\C^*$-equivariant
completion of $V_0$ given by
$$ \bV=\bV_-\cup \bV_+\cup V_0\,,
$$
where
$$
\bV_+=\Spec \cO_\bC[-D_+]\quad \mbox{and}\quad \bV_-=\Spec
\cO_\bC[-D_-]\,.
$$
Moreover, the canonical projection $\pi:V_0\to \bC$ extends to a
$\PP^1$-fibration also denoted $\pi:\bV\to \bC$. The boundary
divisor $\bD=\bV\backslash V_0$ consists of
two disjoint components
$\bC_\pm$ which correspond to the zero sections in $\bV_\pm$,
respectively. \eprop

\bproof We let $\{ p_i \}$ be the points of $\bC$ with
$D_+(p_i)+D_-(p_i)<0$. The fiber over $p_i$ of the orbit map $\pi
: V_0\to \bC$ induced by the inclusion
$\cO_\bC\hookrightarrow\cO_\bC[D_+,D_-]$, is reducible and
consists of
two orbit closures $O_i^\pm$,  see
\cite[I.4]{FlZa1}. Let us consider the $\C^*$-surfaces
$$
V_-=\Spec \cO_\bC[-D_-,D_-]\quad\mbox{and so}\quad  V_+=\Spec
\cO_\bC[D_+,-D_+]\,.
$$
There are natural identifications
$$
V_\pm=V_0\backslash \bigcup_i O_i^\mp \quad\mbox{and}\quad V_+\cup
V_-=V_0\backslash F\,
$$ and open embeddings $V_\pm\hookrightarrow \bV_\pm$,
 where $F$ denotes the fixed point set of the
original $\C^*$-action on $V$. The complements
$\bC_\mp=\bV_\pm\backslash V_\pm$ are the zero sections in
$\bV_\mp$ and so are both isomorphic to $\bC$. Pasting first $V_0$
and $\bV_+$ along their common open subset $V_+$ and gluing
then $\bV_-$ and the resulting surface $V'$ along their common
open subset $V_-$ gives the desired equivariant completion
$\bV$ of $V_0$. \eproof

\brems\label{equivalent} 1. The completion of Proposition
\ref{complet.6} can be constructed with any pair of divisors
$(D_+,D_-)$ on $\bC$. It is not necessary to assume that they are
zero in the points at infinity. For instance, if $p\in
\bC\backslash C$ is a point at infinity and if we replace a pair
of divisors $(D_+,D_-)$ by $(D_+-p, D_-+p)$ then the corresponding
completions $\bV$ and $\bV'$ are easily seen to differ by an
elementary transformation at the fiber $\pi^{-1}(p)\cong \PP^1$.

2. We say that two pairs of $\Q$-divisors
$(D_+,D_-)$ and $(D_+',D_-')$ on $\bC$ are equivalent if
$D_\pm'=D_\pm\pm\div (f)$ for some nonzero meromorphic function
$f$ on $\bC$. By the same arguments as in
\cite[Theorem 4.3(b)]{FlZa1} the hyperbolic $\C^*$-surfaces
$V_0=\Spec
\cO_\bC[D_+,D_-]$ and $V'_0=\Spec O_\bC[D'_+,D'_-]$ over $\bC$ are
equivariantly isomorphic if and only if the pairs $(D_+,D_-)$ and
$(D'_+,D'_-)$ are equivalent on $\bC$. It is easily seen that
equivalent
pairs of divisors on $\bC$ lead to equivariantly isomorphic
completions of $V_0$ in a canonical way.

3. However, starting from equivalent pairs of divisors $(D_+,D_-)$
and $(D_+',D_-')$ on $C$ and extending them by zero in the points
at infinity (as we do in \ref{glue}),  does not lead in
general to equivalent pairs of divisors on $\bC$. Thus the
completion constructed here depends on the equivalence class of
the pair $(D_+,D_-)$ on $C$. Using (1) and (2) it follows that the
completions $\bV$ and $\bV'$ associated to two equivalent pairs
$(D_+,D_-)$ and $(D_+',D_-')$ of $\Q$-divisors on $C$ differ by
elementary transformations at the fibers at infinity.

4.
In the completion $\bV$
the curves $C_\pm$ are parabolic, and $C_+$ is easily seen to be
repulsive whereas $C_-$ is attractive.\erems

Next we describe the singularities of the above completion $\bar
V$ and the intersection pairing on its $\C^*$-invariant divisors.
We use the following notation.

\bsit\label{complet.7} Following \cite[I.4.21]{FlZa1} we let
$\{p_i\}$ be the set of points of $\bC$ with $D_+(p_i)+D_-(p_i)<0$,
and $\{q_j\}$ be the points of $\bC$ with $D_+(q_j)=-D_-(q_j)
\not\in
\Z$. We write
$$D_+(p_i)=-\frac{e_i^+}{m_i^+}\quad\mbox{ and}\quad
D_-(p_i)=\frac{e_i^-}{m_i^-}\quad\mbox{ with}\quad
\gcd(e_i^\pm,m_i^\pm)=1\quad\mbox{ and}\quad\pm m_i^\pm>0\,.
$$
Since $D_+(p_i)+D_-(p_i)<0$ and $m_i^+m_i^-<0$, the
determinant \be\label{dete} \Delta_i=-\left|\ba{ll}e^+_i& e_i^-\\
m_i^+& m_i^-\ea\right| =m_i^+m_i^-(D_+(p_i)+D_-(p_i))\ee is
positive.

The fiber over $p_i$ in $\bV$ consists of two orbit closures
$\bO_i^\pm$ which meet the curves $\bC_\pm$ in points, say,
$p_i^\pm$. Moreover, $\bO_i^+$ and $\bO_i^-$ meet in a unique
point $p_i'$; $p_i^\pm$ and $p_i'$ are the only fixed points over
$p_i$ of the $\C^*$-action on $\bar V$.

Similarly, we let $\bO_j$ be the orbit closure in $\bV$ of the
orbit over $q_j$, and we write
$$D_+(q_j)=-\frac{e_j}{m_j}\quad\mbox{ with}\quad\gcd(e_j,m_j)=1
\quad
\mbox{and}\quad m_j> 0\,.$$ The fiber over $q_j$ is irreducible
and meets $\bC_\pm$ in points $q_j^\pm$.
\esit
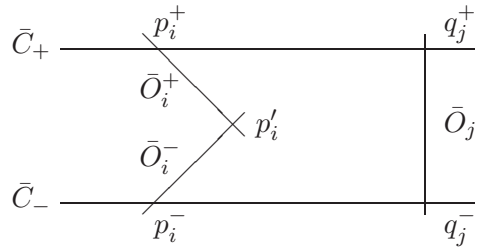
$$
\unitlength0.5cm%
\begin{picture}(12,5.0)
\small \put(1.3,0.8){\line(1,0){11}} \put(1.3,4.9){\line(1,0){11}}
\put(11,0.5){\line(0,1){4.8}}

\put(3.5,5.3){\line(1,-1){2.7}} \put(3.5,0.5){\line(1,1){2.7}}

\put(3.8,5.4){$p_i^+$} \put(3.8,0){$p_i^-$} \put(6.5,2.7){$p_i'$}
\put(3.4,3.6){$\bO_i^+$} \put(3.4,1.8){$\bO_i^-$}

\put(11.5,5.4){$q_j^+$} \put(11.5,0){$q_j^-$}
\put(11.5,2.7){$\bO_j$}

\put(0.0,4.8){$\bC_+$} \put(0.0,0.7){$\bC_-$}
\end{picture}
$$
According to \cite{FlZa1} besides these points $\{p_i,q_j\}$, the
fibers of $\pi$ over all other points are smooth and reduced.

\bsit\label{quosi} Letting further $\Z_d=\langle\zeta\rangle$ be a
cyclic group generated by a primitive $d$-th root of unity
$\zeta$, we consider the $\Z_d$-action on $\A_\C^2$ via
\be\label{la} \zeta.(x,y):=(\zeta x,\zeta^e y)\,, \ee where
$\gcd(d,e)=1$. If $(V,p)\cong (\C^2/\Z_d,\bar 0)$ analytically
then we say that $V$ has a {\em cyclic quotient singularity of
type} $(d,e)$ at $p$. Thus a cyclic quotient singularity of type
$(d,e)$ is also a cyclic quotient singularity of type $(d,\tilde
e)$, where $\tilde e$ is the unique integer with $e\equiv \tilde
e\mod d$ and $0\le \tilde e < d$. The case $d=1$ corresponds to a
smooth point. A cyclic quotient singularity of type $(d,d-1)$ is
an $A_{d-1}$-singularity.\esit

\blem\label{singtype} For an equivariant completion $\bar V$ of
$V_0$ as in Proposition \ref{complet.6}, the following hold.
\begin{enumerate}\item[(a)]
$\bV$ has a cyclic quotient singularity of type $(m_i^+,-e_i^+)$ at
$p_i^+\in\bar V$ and a cyclic quotient singularity of
type $(-m_i^-, -e_i^-)$ at $p_i^-\in\bar V$. In
particular, $p_i^\pm$ is a smooth point of $\bar V$ if and only if
$D_\pm (p_i)$ is integral; that is $m_i^\pm=\pm 1$.
\item[(b)]
$\bV$ has a cyclic quotient singularity of type $(m_j,\mp e_j)$ at
$q_j^\pm\in\bar V$. In particular, $q_j^\pm\in\bar V$
is a  smooth point if and only if $D_\pm (q_j)$ is integral that
is, $m_j=1$. \item[(c)] For a given value of $i$, let $a,b\in\Z$
be integers with
$\left|\ba{ll}a&e_i^+\\
  b&m_i^+\ea\right|=1$, and let
$e^{(i)}=\left|\ba{ll}a&e_i^-\\ b&m_i^-\ea\right|$.
Then $V_0$ has a
cyclic quotient singularity of type $(\Delta_i,e^{(i)})$ at
$p_i'$, see (\ref{dete}). Thus $p_i'\in V_0$ is a smooth point if
and only if $\Delta_i=1$.\end{enumerate} \elem

\bproof The point $p_i^+$ lies on $\Spec \cO_\bC[-D_+]$. As
$-D_+(p_i)=e_i^+/m_i^+$, by \cite[Proposition I.3.8]{FlZa1} $\bV$
has a cyclic quotient singularity of type $(m_i^+,-e_i^+)$ at
$p_i^+$. The other assertions in (a) and (b) follow with the same
argument. For (c) see Theorem 4.15 in \cite[I]{FlZa1}. \eproof

\bprop\label{complet.9} The intersection numbers on $\bV$ are as
follows.

\begin{enumerate}\item[(a)] $\bC_\pm^2=\deg \,D_\pm$ and
$\bC_+.\bC_-=0$.
\item[(b)]
$\bO_i^+.\bC_+=\frac{1}{m_i^+},\quad\bO_i^-.\bC_-=
-\frac{1}{m_i^-}$
and $\,\,\bO_i^+.\bC_-=\bO_i^-.\bC_+=0$.
\item[(c)]
$\bO_j.\bC_\pm=\frac{1}{m_j}$ and $\bO_j^2=0$. \item[(d)]
$\bO_i^+.\bO_i^-=\frac{1}{\Delta_i},\quad
(\bO_i^+)^2=\frac{m_i^-}{\Delta_im_i^+},\quad (\bO_i^-)^2=
\frac{m_i^+}{\Delta_im_i^-}$. \end{enumerate} \eprop

\bproof The first part of (a) follows from Lemma \ref{complet.5},
and the second part is an immediate consequence of the
construction, since $\bC_+$ and $\bC_-$ do not meet.

Again by construction the curves $\bO_i^\pm$ and $\bC_\mp$ do not
meet and so, the intersection numbers $\bO_i^\pm.\bC_\mp$ are equal
$0$. By Proposition 4.18 in \cite[I]{FlZa1} the full fiber over
$p_i$ is given by
$$
\pi^*(p_i)=m^+_i\bO^+_i-m_i^-\bO^-_i\,.
$$
Since its intersection with $\bC_\pm$ is equal to 1, we have
$$ \bO_i^\pm. \bC_\pm=\pm 1/m^\pm_i, $$ proving (b). (c) follows
along the same kind of arguments.

To compute the intersection numbers in (d) we note that the
rational function $ u$ on $\bar V$ as in \ref{complet.0} has a
simple pole along $\bC_+$ and a simple zero along $\bC_-$.
According to Theorem 4.18 in \cite[I]{FlZa1}, the restriction of
$\div \,(u)$ on $V_0$ is given by $\,\,-\sum_j e_j{\bar O}_j-
\sum_i(e^+_i{\bar O}^+_i-e^-_i{\bar O}^-_i)$ and so we obtain
on $\bV$ \be\label{divu} \div\, (u)=-\bC_++\bC_--\sum_j e_j\bO_j-
\sum_i(e^+_i\bO^+_i-e^-_i\bO^-_i). \ee Multiplying with $\bO_i^+$
we get by (b) \be\label{eq1}
e^+_i\bO_i^+\bO^+_i-e^-_i\bO_i^+\bO^-_i=-\bO_i^+\bC_+=-1/m_i^+\,.
\ee As $m_i^+\bO_i^+-m_i^-\bO_i^-$ is numerically equivalent to
any fiber of $\pi$, the product of this divisor with
$m_i^\pm\bO_i^\pm$ is zero. This leads to the equalities
\be\label{eq2} (m_i^+\bO_i^+)^2= (m_i^-\bO_i^-)^2=
(m_i^+\bO_i^+).(m_i^-\bO_i^-). \ee Hence we can rewrite
(\ref{eq1}) as
$$
\Delta_i(\bO_i^+)^2=m_i^-/m_i^+.
$$
Using this and (\ref{eq2}), (d) follows.
\eproof

\subsection{Equivariant resolution of singularities}
In this subsection we consider the minimal resolution of
singularities of $\bV$, which is equivariant by
\ref{equivariant.1}. To describe the boundary divisor we introduce
the following notation.

\bnota\label{complet.10} We abbreviate by a box $\rbh{ \mybox\,}$
with rational weight $e/m$ the weighted linear graph
\be\label{gr1} \cou{C_1}{-k_1}\lin\ldots\lin\cou{C_n}{-k_n}
\qquad=\qquad\boxo{e/m}\ee with $k_1,\ldots,k_n\ge 2$, where $m/
e=[k_1,\ldots, k_n]$, $0< e<m$ and $\gcd(m,e)=1$. A chain of
rational curves $(C_i)$ on a smooth surface with dual graph
(\ref{gr1}) contracts to a cyclic quotient singularity of type
$(m,e)$ \cite{Hi1}. It is convenient to introduce the weighted box
$\quad\boxo{0}\qquad$ for the empty chain. Given extra curves
$E,\,F$ we also abbreviate
\be\label{reversed chain}
\co{E}{}\lin
\co{C_1}\lin\ldots\lin\co{C_n} \qquad
=\qquad\co{E}\llin\boxo{e/m}\qquad = \qquad
\boxo{(e/m)^*}\llin\co{E}\quad\,\, \ee and
\be\label{reversed chain1} \co{C_1}\lin\ldots\lin\co{C_n}
\lin\co{F} \qquad=\qquad
\boxo{e/m}\llin\co{F}\qquad
=\qquad\co{F}\llin\boxo{(e/m)^*}\quad.\ee The orientation of the
chain of curves $(C_i)_i$ in (\ref{gr1}) plays an important role.
Indeed,
 $[k_n,\ldots,k_1]=m/e'$, where $0< e'<m$ and
$ee'\equiv 1 \pmod m$ \cite{Fu, Rus}, and the box $\rbh{
\mybox}\,$ marked with $(e/m)^*:=e'/m$ corresponds to the reversed
chain in (\ref{gr1}).
Note that contracting the curves $C_1,\ldots, C_n$ leads in
both cases to a quotient singularity of type $(m,e)$ on the
ambient surface sitting on $E$ and $F$, respectively,
however with a different orientation; see
e.g.\ \cite[Lemma 5.3.3(1)]{Mi}.
\enota

\bsit\label{complet.11}
Next we consider the minimal resolution of singularities
$\varphi:\tV\to \bV$ of the surface $\bV$. By \cite{OrWa} this
resolution is equivariant and all fibers of
  $\tilde\pi:=\pi\circ \varphi:\tV\to \PP^1$
are chains of rational curves  (cf. also
\ref{equivariant.1}-\ref{equivariant.2}). The proper transforms
$\tC_\pm$ on $\tV$ of the curves $\bC_\pm$ are sections of
$\tilde\pi$. The boundary divisor $\tD=\varphi^{-1}(\bD)$ can be
read off from the following proposition. We recall that
$\{r\}=r-\lfloor r\rfloor$, respectively, $\{D\}=D-\lfloor  D
\rfloor$ stands for the fractional part of a real $r$,
respectively, of a $\Q$-divisor $D$. \esit

\bprop\label{complet.12} \begin{enumerate}[(a)]

\item  The fibers
$F_p=\tilde\pi^{-1}(p)$ in $\tV$ over the points $p\in
\bC\backslash C$ are reduced, isomorphic to $\PP^1$ and satisfy
$F_p.\tC_\pm=1$. Moreover, $F_p.E=0$ for all curves $E$ in $\tD$
different from $\tC_\pm$.

\item $\tC_\pm^2=\deg\, \lfloor
D_\pm\rfloor$.

\item The fibers over the points
$q_j$ together with the curves $\tC_\pm$ are as follows:
$$
\co{\tC_+}\vlin{15}\boxo{\{D_+(q_j)\}}\vlin{15}\co
{\tO_j}\vlin{15}\boxo{\{D_-(q_j)\}^*}\vlin{15}\co{\tC_-}\quad,
$$
where the proper transform $\tO_j$ of $\bO_j$ is a $(-1)$-curve.
All these curves except $\tO_j$ are components of the boundary
divisor $\tD$ of $\tV$.

\item The fibers over the points
$p_i$ together with the curves $\tC_\pm$ are as follows:
$$
\co{\tC_+}\vlin{15}\boxo{\{D_+(p_i)\}}\vlin{15} \co
{\tO_i^+}\lin\co{E_1}\ldots\ldots \co{E_l}\lin\co{\tO_i^-}
\vlin{15}\boxo{\{D_-(p_i)\}^*}\vlin{15}\co{\tC_-}\quad.
$$
Here the chain of rational curves $E_1,\ldots,E_l$ corresponds to
the cyclic quotient singularity at the fixed point $p_i'$ of the
type described in Lemma \ref{singtype}(c). Moreover, the curves
$\tO_i^\pm$ are the proper transforms of $\bO_i^\pm$, and at least
one of them is a $(-1)$-curve. \end{enumerate}\eprop

\bproof (a) is obvious from our construction, see \ref{complet.7}.
By symmetry it is sufficient to prove (b) for the curve $\tilde
C_+$. According to Lemma \ref{complet.9}(a),
$(\bC_+)^2=\deg\,D_+$. By Proposition \ref{complet.2}
$$(\tC_+)^2 =(\bC_+)^2-\sum_{p\in C} \{D_+(p)\}=\deg\,D_+ -
\deg\,\{D_+\}=\deg\,\lfloor D_+(p) \rfloor\,,$$ proving
(b).

To show (c) we may assume that the set $\{p_i,\,q_j\}$ consists of
a single point $q$ so that $D_\pm=\mp e/m[q]$.
By Lemma \ref{singtype}, in this case we deal with the minimal
resolutions of the cyclic quotient singularities of $\bV$ of type
$(m,\mp e)$ at the points $q^\pm\in\bC_\pm$, respectively,
resulting in chains of smooth rational curves with weights
defined via the continuous fraction expansions of $\mp m/e$, see
\ref{complet.10}. As the fiber over $q$ is a chain of
smooth rational
curves, it remains to check that the orientations of
the boxes labeled by $\{D_\pm(q)\}$ are as indicated. By (b),
$\tC_\pm^2=\lfloor D_\pm(q)\rfloor$, and by Proposition
\ref{complet.9}(a), $\bC_\pm^2= D_\pm(q)$. So comparing with
Proposition \ref{complet.2} the orientation of the chain is
indeed as indicated.
Since the fiber $F_q$ can be blown down to a smooth one, one of
its components is a $(-1)$-curve. This can be only the component
$\tO_j$ because the resolution of singularities is minimal.

The proof of (d) is similar and is left to the reader.
\eproof

\brem\label{stable} It is easily seen that any irreducible curve
on $\tV$ stable under the $\C^*$-action on $\tV$ is one of the
curves appearing in the proposition. \erem

\bcor\label{boundarydivisor} If $(\tilde V, \tilde D)$ is the
equivariant  completion of the resolution of singularities of $V$
constructed in \ref{complet.12} and $\tilde\pi : \tV\to \bC$ is
the extension of the orbit morphism $\pi : V\to C$, then the
following hold.
\begin{enumerate} [(a)]\item Every degenerate fiber of the
map $\tilde\pi:\tV\to C$ is a linear chain of rational curves
meeting the sections $\tC_\pm$ in the end components. \item Let
$\bar C$ be a completion of $C$ with ${\rm card}\,(\bC\setminus
C)=s$, and let $\{a_1,\ldots, a_n\}=\{p_i,q_j\}$ be the set of
points of $C$ with $D_+(a_i)\ne 0$ or $D_-(a_i)\ne 0$. Then the
boundary divisor $\tilde D=\tV\setminus V$ has dual graph
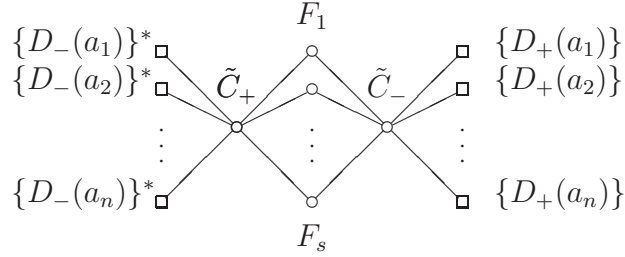
\unitlength1mm
\begin{center}
\begin{picture}(70,35)
      \put(0,25){$\{D_-(a_1)\}^*$}\put(20,25){\mybox}
      \put(0,20){$\{D_-(a_2)\}^*$}\put(20,20){\mybox}
      \put(20,7.5){$\vpoi$}
      \put(0,5){$\{D_-(a_n)\}^*$}\put(20,5){\mybox}
     \put(29.47,15.53){\line(-1,1){8.94}}
     \put(29.3,15.2){\line(-2,1){8.6}}
     \put(29.47,14.47){\line(-1,-1){8.94}}
\put(40,25){$\co{F_1}$} \put(40,20){$\co{}$}
      \put(40,7.5){$\vpoi$}
\put(40,5){$\cou{}{F_s}$}
     \put(49.47,15.53){\line(-1,1){8.94}}
     \put(49.3,15.2){\line(-2,1){8.6}}
     \put(49.47,14.47){\line(-1,-1){8.94}}
\put(30,15){$\co{\tC_-}$}
     \put(30.53,15.53){\line(1,1){8.94}}
     \put(30.7,15.2){\line(2,1){8.6}}
     \put(30.53,14.47){\line(1,-1){8.94}}
\put(30,15){$\co{\tC_+}$}

\put(50,15){$\co{\tC_-}$}
     \put(50.53,15.53){\line(1,1){8.94}}
     \put(50.7,15.2){\line(2,1){8.6}}
     \put(50.53,14.47){\line(1,-1){8.94}}
      \put(60,25){$\mybox\quad \{D_+(a_1)\} $}
      \put(60,20){$\mybox\quad  \{D_+(a_2)\} $}
      \put(60,7.5){$\vpoi$}
      \put(60,5){$\mybox\quad \{D_+(a_n)\}$}
\end{picture}
\end{center}
\vskip0.3truecm Besides possibly $\tC_\pm\cong \bC$ all the curves
are rational, and $F_1,\ldots,F_s$ are the fibers over the points
at infinity.

\item In particular, if $C=\A^1$ then the boundary divisor $\tilde
D$ consists of smooth rational curves, and the dual graph $\G
(\tilde D)$ is \unitlength1mm
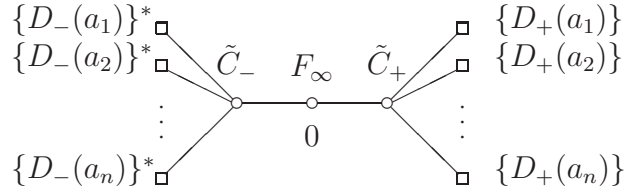
\begin{center}
\begin{picture}(70,35)
      \put(0,25){$\{D_-(a_1)\}^*$}\put(20,25){\mybox}
      \put(0,20){$\{D_-(a_2)\}^*$}\put(20,20){\mybox}
      \put(20,7.5){$\vpoi$}
      \put(0,5){$\{D_-(a_n)\}^*$}\put(20,5){\mybox}
     \put(29.47,15.53){\line(-1,1){8.94}}
     \put(29.3,15.2){\line(-2,1){8.6}}
     \put(29.47,14.47){\line(-1,-1){8.94}}
      \put(60,25){$\mybox\quad \{D_+(a_1)\} $}
      \put(60,20){$\mybox\quad  \{D_+(a_2)\} $}
\put(30,15){$\co{\tC_-}\lin\cou{F_\infty}{0}\lin\co{\tC_+}$}
      \put(60,7.5){$\vpoi$}
      \put(60,5){$\mybox\quad \{D_+(a_n)\}$}
     \put(50.53,15.53){\line(1,1){8.94}}
     \put(50.7,15.2){\line(2,1){8.6}}
     \put(50.53,14.47){\line(1,-1){8.94}}
\end{picture}
\end{center}

\item $V$ is a Gizatullin surface, i.e.\ the dual graph $\G
(\tilde D)$ is a linear chain of rational curves, if and only if
$C\cong\A^1_\C$ and each of the fractional parts $\{D_\pm\}$ is
either zero or supported at one point: $\supp (\{D_\pm\})
\subseteq \{p_\pm\}$.

\item The dual graph of $\tD$ is circular if and only if $s=2$ and
the divisors $D_\pm$ are integral. Moreover, in this case the dual
graph is $((0, \deg D_+,0,\deg D_-))$, which has standard form
$((0,0,0,\deg\,(D_++D_-)\,))$.
\end{enumerate}\ecor

\brem\label{hircyc} We note that every surface as in (e) can be
obtained from a Hirzebruch surface by blowing up at some distinct
points of two disjoint sections (not at the same fiber) and
deleting two other fibers and the proper transforms of these
sections. The $\C^*$-action is vertical and the sections are
parabolic curves. \erem

\bexas\label{complet.13} 1. It can happen that both $\tO_i^\pm$
are $(-1)$-curves. Indeed, assume that for some $i$ the
coefficients $D_\pm(p_i)$ at $p_i$ are both integral and
$-n:=D_+(p_i)+D_-(p_i)< 0$. In this case by Lemma \ref{singtype}
the points $p_i^\pm\in \bC_\pm$ are smooth and $p_i'\in\bV$ is a
cyclic quotient singularity of type $(n, n-1)$ (with
$\Delta_i=n$). Since $n/(n-1)=[2,\ldots,2]$ ($n-1$ times) the
fiber of $\tV\to \bC$ over $p_i$ together with the curves
$\tC_\pm$ is
$$ \cou{\tC_+}{}\lin\cou{\tO_i^+}{-1}\llin\boxo{A_{n-1}}
\llin\cou{\tO_i^-}{-1}\lin\cou{\tC_-}{} $$ with a chain
$A_{n-1}=[[(-2)_{n-1}]]$ of $(-2)$-curves of length $n-1$ in the
middle. Indeed, by Proposition \ref{complet.9}(d),
$(\bO_i^\pm)^2=-1/n$. Applying Proposition \ref{complet.2} we
obtain $(\tO_i^\pm)^2=\lfloor -1/n \rfloor=-1$.

2. Let $C$ be a nodal cubic in $\PP^2$. We claim that the smooth
affine surface $V=\PP^2\setminus C$ does not admit a
$\C^*$-action. Indeed, $C$ has dual graph $((9))$ with standard
form $((0,0,(-2)_6,-3))$, and so this graph is not birationally
equivalent to a one in (e) above. Henceforth $V$ does not admit a
hyperbolic $\C^*$-action. We will see below that the dual graphs
of equivariant completions of parabolic and elliptic
$\C^*$-surfaces are trees, which excludes existence of a parabolic
or elliptic $\C^*$-action on $V$. \eexas

\subsection{Parabolic and elliptic $\C^*$-surfaces}

In this subsection we give a short description of the boundary
divisors of parabolic and elliptic $\C^*$-surfaces.

\bsit\label{paracase} {\it Parabolic case.} We let $V=\Spec
A_0[D]$ be a parabolic $\C^*$-surface, where $D$ is a $\Q$-divisor
on a smooth affine curve $C=\Spec A_0$. The projection $A_0[D]\to
A_0$ provides a section $\iota: C\to V$ with image $C_0=\iota (
C)$.

We recall that $V$ has a cyclic quotient singularity of type
$(m,e)$ at $\iota (p)\in C_0$ if $D(p)=-e/m$, see \cite[Prop.
I.3.8]{FlZa1}.

Letting as before $\bC$ be a smooth completion of $C$ with $s$
points at infinity, we consider $D$ as a $\Q$-divisor on $\bC$ and
we identify the function field $K=\Frac(C)$ with the constant
sheaf $K$ on $\bC$. We form a sheaf of $\cO_\bC$-algebras
$$
\cO_\bC[D]\subseteq K[u,u^{-1}]
$$
as in \ref{complet.5}. The
corresponding normal $\C^*$-surface $V_0=\Spec \cO_\bC[D]$ can be
completed as follows.\esit

\bprop\label{pcomplet.6} There is a natural $\C^*$-equivariant
completion of $V$ given by
$$
\bV=V_0\cup V_\infty,
$$
where $V_0$ and $V_\infty=\Spec \cO_\bC[-D]$ are pasted along
$V^*:= V_0\cap V_\infty =\Spec \cO_\bC[D,-D]$ via $u\mapsto
u^{-1}$. Moreover, the canonical projections $\pi:V_0\to \bC$ and
$\pi:V_\infty\to\bC$ coincide on the intersection and so provide a
$\PP^1$-fibration also denoted $\pi:\bV\to \bC$. The  boundary
divisor $\bD=\bV\backslash V$ has a decomposition
$$
\bD=\bC_\infty\cup F_1\cup\ldots\cup F_s,
$$
where $F_1,\ldots, F_s$ are the fibers of $\pi$ over $\bar
C\backslash C$ and $\bC_\infty$ corresponds to the section in
$\bV_\infty$ induced by the projection  $\cO_\bC[-D]\to\cO_\bC$.
\eprop

\bproof By Proposition 4.1 and Remark 4.20 in \cite[I]{FlZa1},
$V^*=\Spec \cO_\bC[D,-D]$ can be identified with the open subset
$V_0\backslash \bC_0$ of $V_0$, and similarly, with the open
subset $V_\infty\backslash \bC_\infty$ of $V_\infty$. Thus pasting
$V_0$ and $V_\infty$ along $V^*$ gives an equivariant completion
of $V$, cf. Proposition \ref{complet.6}. The above description of
the boundary divisor is now straightforward. \eproof

We let further $\tV$ be the minimal resolution of singularities of
$\bV$, and  $\tC_0$, $\tC_\infty$ be the proper transforms of the
sections $\bC_0$ and $\bC_\infty$, respectively. For every point $
p\in C$ with $D(p)=-e/m$ the surface $\bV$ has a cyclic quotient
singularity of type $(m,m-e)$ at the point $p'\in\tC_\infty$ over
$p$, cf. Lemma \ref{singtype}(a). Thus using
\ref{boundarydivisor}(b) the dual graph of the boundary divisor is
as follows: \unitlength1mm
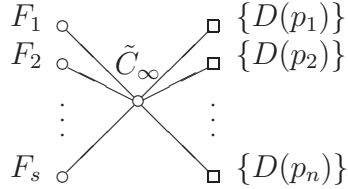
\begin{center}
\begin{picture}(60,31)
  \put(43,25){$\{D(p_1)\}$}\put(40,25){\mybox}
     \put(43,20){$\{D(p_2)\}$}\put(40,20){\mybox}
      \put(40,7.5){$\vpoi$}
    \put(43,5){$\{D(p_n)\}$}\put(40,5){\mybox}
    \put(29.47,15.53){\line(-1,1){8.94}}
    \put(29.3,15.2){\line(-2,1){8.6}}
    \put(29.47,14.47){\line(-1,-1){8.94}}
   \put(20,25){\co{}} \put(13,25){$F_1$}
\put(20,20){\co{}}   \put(13,20){$F_2$}
  \put(20,7.5){$\vpoi$}
 \put(20,5){\co{}}  \put(13,5){$F_s$}

    \put(30.53,15.53){\line(1,1){8.94}}
    \put(30.7,15.2){\line(2,1){8.6}}
    \put(30.53,14.47){\line(1,-1){8.94}}
\put(30,15){$\co{\tC_\infty}$}
\end{picture}
\end{center}
where $\{p_i\}$ are the points of $C$ with $\{D(p_i)\}\neq 0$.
Thus the dual graph of the boundary divisor $D=\tV\setminus V$ is
a linear chain of rational curves if and only if $C\cong \A^1_\C$
and ${\rm supp}\,(\{D\})$ is either empty or consists of
one point.

\bsit\label{ellcase} {\it Elliptic case.} We let now $V=\Spec A$,
where $A=\bigoplus_{i\ge 0} A_i$ with $A_0=\C$ is a positively
graded normal 2-dimensional $\C$-algebra of finite type. So $V$ is
an elliptic $\C^*$-surface. By the results of Dolgachev and
Pinkham, see \cite[I]{FlZa1}, the projective curve $C=\Proj A$ is
smooth, and there is a $\Q$-divisor $D$ on $C$ with $\deg D >0$
such that
$$A_n=H^0\left(C, \cO_C\lfloor nD\rfloor \right)\cdot u^n\subseteq
\Frac (\cO_C)[u],\quad \forall n\ge 0\,.$$ The elliptic
$\C^*$-surface $V$ can be obtained in the following way. Consider
the surface
$$S_0=\Spec (\cO_C[D])$$ with a parabolic
$\C^*$-action provided by the grading of $\cO_C[D]$. The inclusion
$\cO_C\hookrightarrow \cO_C[D]$ gives the orbit map $S_0\to C$,
and the projection $\cO_C[D]\to\cO_C$ gives a section $\iota :
C\to S_0$. The natural map $A\to \cO_C[D]$ yields a morphism $\pi:
S_0\to V$, which is the contraction of the curve $C_0=\iota (C)
\hookrightarrow S_0$. As in the parabolic case, $S_0$ has a cyclic
quotient singularity of type $(m,e)$ at $\iota (p) \in C_0$, where
$D(p)=-e/m$. We obtain now a completion $\bar {S_0}$ of $S_0$ as
follows.\esit

\bprop\label{ellca} There is a natural $\C^*$-equivariant
completion $\bar S$ of $S_0$ given by $$\bar S=S_0\cup
S_{\infty}\,,$$ where $S_0$ and $S_{\infty}=\Spec \cO_C[-D]$ are
pasted along $S_0\cap S_{\infty}=\Spec \cO_C[D,-D]$ via $u\mapsto
u^{-1}$. The canonical projections $S_0\to C$ and $S_{\infty}\to
C$ provide a projection $\pi: \bar S\to C$, and the section
$C_{\infty}=\bar S\setminus S_0\subseteq S_{\infty}$ of $\pi$ is
induced by the projection $\cO_C[-D]\to\cO_C$. \eprop

The proof is the same as in the parabolic case. Consider further
the minimal resolution of singularities $\sigma : \tilde S\to \bar
S$, and let $\tilde C_{\infty}$ be the proper transform of
$C_{\infty}$. For every point $p\in C$ with $D(p)=-e/m$ the
surface $\bar S$ has a cyclic quotient singularity of type $(m,
m-e)$ at the point $p'\in C_{\infty}$ over $p$. Thus similarly as
before the boundary divisor $\tilde S\setminus S_0$ has dual graph
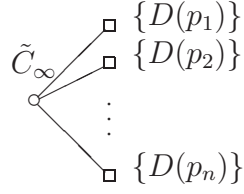
\unitlength1mm
\begin{center}
\begin{picture}(60,31)
  \put(43,25){$\{D(p_1)\}$}\put(40,25){\mybox}
     \put(43,20){$\{D(p_2)\}$}\put(40,20){\mybox}
      \put(40,7.5){$\vpoi$}
    \put(43,5){$\{D(p_n)\}$}\put(40,5){\mybox}


 \put(30.53,15.53){\line(1,1){8.94}}
    \put(30.7,15.2){\line(2,1){8.6}}
    \put(30.53,14.47){\line(1,-1){8.94}}
\put(30,15){$\co{\tilde C_\infty}$}
\end{picture}
\end{center} where $(p_i)$ are the points of $C$ with
$\{D(p_i)\}\neq 0$.

Since $V$ is obtained from $S_0$ by contracting $C_0$,
contracting $C_0$ on $\bar S$ yields a completion $\bar V$ of $V$.
The minimal resolution of singularities $\tilde V\to \bar V$ of
$\bar V$ is also equivariant, and the boundary divisor
$\tV\setminus V$ is as shown in the above diagram. This
divisor is a linear chain of rational curves provided that $C$ is
rational and $\{D\}$ is concentrated in at most 2 points.

  \section{Boundary zigzags of
Gizatullin $\C^*$-surfaces}

In this section we address Gizatullin surfaces. By definition
(see the Introduction) these are normal affine surfaces admitting
completion by a zigzag, i.e.\ by an SNC divisor whose components
are rational curves and the dual graph is linear.

\subsection{Smooth Gizatullin surfaces}
By Theorem 2.15 in \cite{FKZ} any Gizatullin surface admits a
completion with a standard zigzag $[[0,0,w_2,\ldots,w_n]]$, $n\ge
1$, as boundary: \be\label{zigzag}\,\qquad
\cou{C_0}{0}\lin\cou{C_1}{0}\lin \cou{C_2}{w_2}\lin\ldots\lin
\cou{C_n}{w_n}\quad, \qquad \mbox{where}\qquad\begin{cases} w_i\le
-2\,\,\,\,\,\forall i\,\,\,&\mbox{ if}\,\,\, n\ge 3\\ w_2\le
0,\,\,w_2\neq -1\,&\mbox{ if}\,\,\, n=2\,.\,\end{cases}\ee By
Corollary 3.5 in \cite{FKZ} this zigzag is unique up to reversing
the sequence of weights  $(w_2,\ldots, w_n)$. The following lemma
shows that actually every such zigzag can be the boundary of a
smooth Gizatullin surface.

\blem\label{boundary} (\cite[I]{Gi} or also \cite[I]{Du2}) Every
standard zigzag (\ref{zigzag}) occurs as boundary divisor of a
smooth Gizatullin surface $X$. \elem

\bproof We start with the quadric $Q=\PP^1\times\PP^1$ and the
curve $C_0+C_1+C_2$ on $Q$, where \be\label{curves3}
C_0=\{\infty\}\times\PP^1,\quad C_1=\PP^1\times\{\infty\}\quad
\mbox{and}\quad C_2=\{0\}\times\PP^1\,.\ee In case $n=1$ we let
$X=Q$ and $D=C_0+C_1$ with $C_0,C_1$ as above. If $n\ge 2$ then
performing a sequence of outer blowups over a point $x_0\in
C_2\backslash C_1$ we obtain a linear chain of rational curves
$D=C_0+C_1+\ldots+C_n$ with dual graph $Z=[[0,0,0]]$ if $n=2$
(here no blowup is necessary) and $Z=[[0,0,-1,(-2)_{n-3},-1]]$ if
$n\ge 3$, respectively. Performing further blowups with centers at
distinct points of the curves $C_i$ different from the double
points of $D$, we can achieve the prescribed weights $C_i^2=w_i\le
-2$, $i=2,\ldots,n$.

Letting $\bar X$ be the resulting smooth projective surface
dominating $Q$, we denote by $\bar D= \bar C_0+\bar
C_1+\ldots+\bar C_n$ the proper transform of $D$ in $\bar X$. It
remains to check that the smooth open surface $X=\bar X\backslash
\bar D$ is affine. For this it is enough to show that, for a
sequence of positive multiplicities $m_0,\ldots,m_n$, the divisor
$D'=\sum_{i=0}^n m_i\bar C_i$ on $\bar X$ is ample. It is easily
seen that $\bar D$ meets every irreducible curve $C$ on $\bar X$
different from all the $\bar C_i$, hence $C\cdot \bar D'>0$. Also
$\bar C_i\cdot D'>0$ for every $i=0,\ldots,n$ provided that
$m_{i+1}+m_{i-1}>-m_iw_i$ for all $i$. The latter can be achieved
recursively starting with $m_n=1$. Now such a divisor $D'$ is
ample by the Nakai-Moishezon criterion. \eproof

\subsection{Toric Gizatullin surfaces}
In this part we answer the question what further restrictions on
the boundary zigzag of a Gizatullin surface are imposed by the
presence of a $\C^*$-action. The answer provided by Proposition
\ref{prpzigzag} and Theorem \ref{thmzigzag} below depends on
whether the surface is smooth or not. Let us first examine the
toric case.

\blem\label{torics} \begin{enumerate}[(a)] \item Every smooth
toric affine surface is isomorphic either to $\C^*\times\C^*$, to
$\A^1_\C\times \C^*$ or to $\A^2_\C$. Every normal singular toric
affine surface is isomorphic to $V_{d,e}:=\A^2/\langle \zeta
\rangle$, where the primitive $d$-th root of unity $\zeta$ acts on
$\A^2$ via $\zeta . (x,y)=(\zeta x, \,\zeta^e y)$ for some
$d>1,\,e\in \Z$ with $\gcd (e,d)=1$.

\item A $\C^*$-surface $V=\Spec A_0[D_+,D_-]$ with
$A_0=\C[t]$ is toric if and only if $(D_+,D_-)\sim
\left(-\frac{e^+}{m^+} [p_0], \, \frac{e^-}{m^-} [p_0]\right)$ for
some point $p_0\in\A^1$.
\end{enumerate}\elem

\bproof (a) is well known and can be found in e.g.\ \cite[I,
Example 2.3 and II, Example 2.8]{FlZa1}. To deduce (b), if
$V=\Spec A_0[D_+,D_-]$ and $(D_+,D_-)\sim \left(-\frac{e^+}{m^+}
[p_0], \, \frac{e^-}{m^-} [p_0]\right)$, then $V$ is toric as was
shown in the proof of Theorem 4.15(c) in \cite[I]{FlZa1}.
Conversely assume that for some pair of $\Q$-divisors $(D_+,D_-)$
the surface $V=\Spec A_0[D_+,D_-]$ is toric. According to
\cite{FlZa2}, Theorem 4.5 and its proof the pair $(D_+,D_-)$ has
then the claimed form, so the lemma follows. \eproof

\bprop\label{prpzigzag}
\begin{enumerate}[(a)]
\item Any standard zigzag $(\ref{zigzag})$
occurs as the boundary divisor of a normal toric affine
surface\footnote{Hence also of a normal surface with a hyperbolic
(elliptic, parabolic) $\C^*$-action.}.

\item  A standard zigzag $(\ref{zigzag})$ occurs as the boundary
divisor of a smooth toric affine surface if and only if it is
$[[0,0]]$ or $[[0,0,0]]$.
\end{enumerate}\eprop

\bproof To show (a), given a standard zigzag
$[[0,0,w_2,\ldots,w_n]]$ as in (\ref{zigzag}) we write
$$\frac{m}{e}=[-w_2+1, -w_3,\ldots, -w_n]\quad\mbox{ with}\quad
\gcd(e,m)=1\,.$$ We also consider the pair of $\Q$-divisors
$(D_+,D_-)=(-\frac{e}{m}[0],0)$ on the affine line $C=\A^1$. By
Lemma \ref{torics}(b), $V=\Spec A_0[D_+,D_-]$ with $A_0=\C[t]$ is
a toric surface. According to Corollaries \ref{equivariant.5},
\ref{boundarydivisor}(c) and Proposition \ref{complet.12}(d), $V$
has a $\T$-equivariant completion $\tV$ with boundary divisor
$$\tD'\,\,=\quad
\cou{\tilde C_-}{0}\lin\cou{F_\infty}{0}\lin\cou{\tilde C_+}{-1}
\vlin{15} \boxo{\frac{e}{m}}\quad\quad =\quad\quad \cou{\tilde
C_-}{0}\lin\cou{F_\infty}{0}\lin\cou{\tilde C_+}{-1} \llin
\cou{}{w_2-1}\llin\cou{}{w_3}\lin\ldots\lin\cou{}{w_n}\quad.
$$
Contracting $\tilde C_+$ we perform a $\T$-equivariant outer
elementary transformation which consists in blowing up at the only
fixed point on $\tilde C_-\ominus F_\infty$ of the torus action on
$\tV$ and then blowing down the proper transform of $\tilde C_-$.
This results in a new equivariant completion of $V$ with the given
standard zigzag $[[0,0,w_2,\ldots,w_n]]$ as boundary, proving (a).

Now (b) follows from Lemma \ref{torics}(a) by virtue of the
uniqueness (up to reversion) of a standard zigzag in its
birational equivalence class, see Corollary 3.5 in \cite{FKZ}.
\eproof

\subsection{Smooth Gizatullin $\C^*$-surfaces}
\bthm\label{thmzigzag} A standard zigzag occurs as the boundary
divisor of a smooth affine hyperbolic $\C^*$-surface if and only
if it can be written in one of the forms $[[0,0]],\,\,[[0,0,0]]$,
$$
\cou{}{0}\lin\cou{}{0}\llin
\boxo{(e_1/m_1)^*}\llin\cou{}{-2-k}\llin\boxo{e_2/m_2} \qquad\quad
\mbox{or }\quad(ii)\qquad
\cou{}{0}\lin\cou{}{0}\lin\boxo{A_{m_1-1}}\lin\cou{}{-2-k}
\lin\boxo{A_{m_2-1}} \quad,\leqno (i)$$ where as before
$A_k=[[(-2)_k]]$, $k\ge 0$, $m_i\ge 1$, $\gcd(e_i,m_i)=1$ for
$i=1,2$, and either \be\label{toon}
\frac{e_1}{m_1}+\frac{e_2}{m_2}=1
  \quad\mbox{ or
}\quad \frac{e_1}{m_1}+\frac{e_2}{m_2}=1-\frac{1}{m_1m_2}\,. \ee
\ethm

\bproof We suppose first that $V=\Spec A_0[D_+,D_-]$ is a smooth
affine surface with a hyperbolic  $\C^*$-action, completed by a
standard zigzag. By Corollary \ref{boundarydivisor}(d) $A_0\cong
\C[t]$ and the support of each of the fractional parts $\{D_\pm\}$
is  empty or consists of just one point $p_\pm$. Actually we
establish below that (i) holds if $p_+=p_-$ and (ii) holds if
$p_+\neq p_-$.

If $V$ is a smooth toric surface then  the assertion follows from
Proposition \ref{prpzigzag}(b). So we may assume for the rest of
the proof that $V$ is not toric.

\no$\bullet$  Suppose first that the fractional parts $\{D_\pm\}$
are supported at  the same point $p_+=p_-$ or that one or both of
them are zero. By a coordinate change of the base and passing to
an equivalent pair $(D_+,D_-)$ we may assume that
$p_\pm=0\in\A^1_\C$ and
$$\left(D_+,D_-\right)=\left(\left(\frac{e_1}{m_1}-1\right)[0],\,
\frac{e_2}{m_2}[0]-D'\right)\,\,\,\mbox{with}\,\,\, 0\le
e_2<m_2,\,\,\, \gcd(e_i,m_i)=1,\,\,i=1,2\,,$$ where $D'$ is an
effective integral divisor of degree, say, $ k+1\ge 0$ with
$0\not\in\supp (D')$. Actually $k\ge 0$ since otherwise, $D_\pm$
being concentrated at one point,  by Lemma \ref{torics}(b) $V$
would be a smooth toric surface, which has been excluded.

The fibers of $\pi : V\to\A^1_\C$ over the points $p_i\in {\rm
supp}\,D'$ are reducible and singular at the points $p_i'$, see
\ref{complet.7}. According to Lemma \ref{singtype}(c) $p_i'\in V$
is a smooth point if and only if $\Delta_i=D'(p_i)=1$. Since $V$
is supposed to be smooth, $D'$ is supported at $ k+1$ distinct
points.

\noindent $\diamond$ The fiber in $V$ over $0\in\A^1_\C$ is
irreducible (and so $V$ is automatically smooth along this fiber)
if and only if
$$D_+(0)+D_-(0)=0\quad\Longleftrightarrow\quad \frac{e_1}{m_1}+
\frac{e_2}{m_2}=1\,.$$ The latter agrees with the first equality
in (\ref{toon}). Note that this is also true if $m_1=1$ or $m_2=1$
since in this case the corresponding boxes in (i) are empty.

\noindent $\diamond$ The fiber over $0$ is reducible if and only
if $D_+(0)+D_-(0)<0$ i.e., $e_1/m_1+e_2/m_2<1$. Moreover
$$D_+(0)=\frac{e_1-m_1}{m_1}=-\frac{e_0^+}{m_0^+}
\quad\mbox{and}\quad
D_-(0)=\frac{e_2}{m_2}=\frac{e_0^-}{m_0^-}\,,$$ so again by Lemma
\ref{singtype}(c), $V$ is smooth along this fiber if and only if
$$\Delta_0=-\left| \begin{array}{cc}
  e_0^+ & e_0^- \\
  m_0^+ & m_0^- \\
\end{array}\right|=-\left| \begin{array}{cc}
  m_1-e_1 & -e_2 \\
  m_1 & -m_2 \\
\end{array}\right|=1\quad\Longleftrightarrow\quad
\frac{e_1}{m_1}+ \frac{e_2}{m_2}=1-\frac{1}{m_1m_2}\,.$$ The
latter agrees with the second equality in (\ref{toon}).

By Proposition \ref{complet.12}(b) we have \be\label{squares}
\tC_+^2=\deg\,\lfloor D_+\rfloor=-1\quad\mbox{ and}\quad
\tC_-^2=\deg\,\lfloor D_-\rfloor=-1-k\,.\ee Thus by virtue of
Corollary \ref{boundarydivisor} the boundary divisor of the
completion $\tV$ constructed in \ref{complet.12} has the form
$$ \boxo{(\frac{e_1}{m_1})^*}\vlin{15}\cou{\tC_+}{-1} \lin
\cou{F_\infty}{0}
\llin\cou{\tC_-}{-1-k}\vlin{15}\boxo{\frac{e_2}{m_2}}\qquad.
$$ Performing an
elementary transformation at $F_\infty\cap \tC_-$ by blowing up
this point and blowing down the proper transform of $F_\infty$, we
arrive at a linear chain with two zero weights in the middle. By
virtue of Lemma 2.12(a) in \cite{FKZ} applying further a sequence
of elementary transformations we can move this pair of zero
weights to the left to obtain a standard zigzag of type (i).

\noindent $\bullet$ If now $\{D_\pm\}\neq 0$ and  $p_+\neq p_-$
then we can write
$$D_+=\frac{e_1}{m_1}[p_+]\quad\mbox{ and}\quad D_-=
\frac{e_2}{m_2}[p_-]-D',\qquad\mbox{ where} \quad m_1,\,m_2\ge 2
$$
and $D'$ is an effective integral divisor of degree $k\ge 0$,
whose support does not contain the points $p_\pm$. As before, the
condition that $V$ is smooth forces by Lemma \ref{singtype}(c)
that $D'$ is supported at $k$ distinct points and $e_1=e_2=-1$, so
that
$$D_+=\frac{-1}{m_1}[p_+]\quad\mbox{ and}\quad D_-=
\frac{-1}{m_2}[p_-]-\sum_{i=1}^k p_i\quad\mbox{with}\quad p_i\neq
p_\pm\,\,\forall i\,.$$ Again (\ref{squares}) hold and so, the
boundary divisor $\tD$ of the smooth equivariant completion $\tV$
of $V$ is in this case
$$
\boxo{(\frac{m_1-1}{m_1})^*}\vlin{18}\cou{\tC_+}{-1} \lin
\cou{F_\infty}{0}
\llin\cou{\tC_-}{-1-k}\vlin{18}\boxo{\frac{m_2-1}{m_2}}\quad,
$$ see Examples \ref{complet.3} and \ref{complet.13}.
Performing a sequence of inner elementary transformations we can
transform this into the standard zigzag (ii), as required.

\noindent $\bullet$ Vice versa, given a linear chain $\G$ as in
(i) or (ii), we choose the divisors $D_\pm$ as in the proof above.
This yields a smooth affine surface $V=\Spec A_0[D_+,D_-]$ with
$A_0=\C[t]$ equipped with a hyperbolic $\C^*$-action, which admits
an equivariant completion by a standard zigzag $\tD_{\rm st}$ with
dual graph $\G$. \eproof

\brem\label{switch} Reversing the grading on $A_0[D_+,D_-]$ or,
equivalently, switching $\lambda\longmapsto\lambda^{-1}$ in the
$\C^*$-action amounts to interchanging $D_+$ and $D_-$. This also
amounts to reversing the standard zigzags in (i) or in (ii).\erem

The following corollary is similar to Russell's description of the
Ramanujam-Morrow graphs \cite[3.3]{Rus}.

\bcor\label{corzigzag} A standard zigzag $[[0,0,w_2,\ldots,w_n]]$
occurs as the boundary divisor of a smooth $\C^*$-surface if and
only if one of the following conditions is satisfied.
\begin{enumerate}
\item[(i$\,'$)] For some $i$ with $2\le i\le n$, the zigzag  $[[
w_2,\ldots,w_{i-1},-1,w_{i+1}\ldots,w_n]]$ is contractible to
$[[0]]$ or to  $[[ -1]]$. \item[(ii$\,'$)] $[[ 0,0,
w_2,\ldots,w_n]] = [[ 0,0,(-2)_{\alpha},\,-2-k,\,(-2)_{\beta}]]$
for some $\alpha,\beta,k\ge 0$ with $\alpha+\beta=n-2$.
\end{enumerate}
\ecor

\bproof We must show that (i$'$) is equivalent to condition
(i) of
Theorem \ref{thmzigzag}(b). Consider first the case that
$e_1/m_1+e_2/m_2=1$ in \ref{thmzigzag}(b)(i). This means that
$m:=m_1=m_2$ and $e_2=m-e$, where $e:=e_1$. Replacing in the
zigzag from \ref{thmzigzag}(b)(i) the weight $-2-k$ by $-1$ and
choosing $e'$, $0\le e'< m$, with $ee'\equiv 1\mod m$, by virtue
of \ref{complet.10} we obtain
$$
\boxo{(\frac{e}{m})^*}\vlin{15}\co{-1}\vlin{15}
\boxo{\frac{m-e}{m}}
\qquad\qquad\mbox{=}\qquad\qquad
\boxo{\frac{e'}{m}}\vlin{15}\co{-1}\vlin{18}
\boxo{(\frac{m-e'}{m})^*}\quad.
$$ Letting now in Proposition
\ref{complet.12}(c) $D_+=e'/m[0]$ and $D_-=-e'/m[0]$, the latter
chain occurs as the dual graph of the fiber over $0\in\A^1$ of a
$\PP^1$-fibration $\tilde \pi:\tV\to\PP^1$ on a smooth surface
$\tV$. Therefore it contracts to $[[0]]$.

Similarly, by Proposition 4.9(b) in \cite{FKZ}, the condition
\be\label{eqfrac}
\frac{e_1}{m_1}+\frac{e_2}{m_2}=1-\frac{1}{m_1m_2} \ee is
equivalent to the contractibility to $[[ -1]]$ of the graph
$$\boxo{(\frac{e_1}{m_1})^*}\vlin{15}\co{-1}\vlin{15}
\boxo{\frac{e_2}{m_2}} \quad.$$ Now the proof is completed.
\eproof

\bexas\label{eboundary.3} 1. Every zigzag $[[0,0,w_2]]$ with
$w_2\le -2$ is of type (ii) in Theorem \ref{thmzigzag}(b) with
$m_1=m_2=1$, so that the boxes labelled by $A_{m_i-1}$, $i=1,2$,
are empty.

2. The zigzag $[[0,0,w_2,w_3]]$ with $w_2,w_3\le -2$
satisfies one
of the conditions in Corollary \ref{corzigzag} (and so,
corresponds to a smooth Gizatullin surface with a hyperbolic
$\C^*$-action) if and only if at least one of the weights
$w_2,w_3$ is equal to $-2$.

Indeed, the linear chains $[[-1,w_3]]$ and $[[w_2,-1]]$ cannot be
contracted to $[[0]]$ whatever are the weights $w_i\le -2$,
$i=1,2$. Moreover, under the above condition, and only then, one
of these chains contracts to $[[-1]]$ and so, the zigzag
$[[0,0,w_2,w_3]]$ satisfies (i$'$) and, simultaneously, (ii$'$).

3. A graph $[[0,0,w_2,w_3,w_4]]$ ($w_i\le -2$) corresponds to a
smooth Gizatullin surface with a hyperbolic $\C^*$-action if and
only if either two of the weights $w_2,w_3,w_4$ are equal to $-2$
or $(w_2,w_4)$ is one of the pairs $(-2,-3)$ or $(-3,-2)$.

Indeed, in the first case (ii$'$) in Corollary \ref{corzigzag} is
fulfilled, and in the second one (i$'$) holds. Actually the chain
$[[w_2,-1,w_4]]$ contracts to $[[-1]]$ (to $[[0]]$, respectively)
if and only if $(w_2,w_4)$ is one of the pairs $(-2,-3)$ or
$(-3,-2)$ ($(-2,-2)$, respectively). Moreover, the chains
$[[-1,w_3,w_4]]$ and $[[w_2,w_3,-1]]$ cannot be contracted to
$[[0]]$, and they are contracted to $[[-1]]$ if and only if
$w_3=w_4=-2$, respectively, $w_2=w_3=-2$.\eexas

An elliptic or parabolic Gizatullin $\C^*$-surface is necessarily
toric, see Corollary 4.4 in \cite[II]{FlZa1}. In particular, if
such a surface is smooth then it is equivariantly isomorphic to
$\A^2$ or $\A^1\times\C^*$ with a linear $\C^*$-action. Therefore
the above examples and Lemma \ref{boundary} lead to the following
corollary.

\bcor\label{noaction} There exist smooth Gizatullin surfaces that
do not admit any $\C^*$-action.\ecor

\brem Every Gizatullin surface admits two non-conjugated
$\C_+$-actions. However \cite[Corollary 3.4]{FlZa2} if a normal
affine surface $V\not\cong \C^*\times\C^*$ admits two distinct, up
to switching $\lambda\longmapsto \lambda^{-1}$ in one of them,
$\C^*$-actions then it also admits a $\C_+$-action. Moreover by
\cite[Theorem 3.3]{FlZa2} $V$ is a Gizatullin surface provided
that these $\C^*$-actions are non-conjugate and remain
non-conjugate after switching $\lambda\longmapsto \lambda^{-1}$ in
one of them.\erem

\section{Extended graphs of Gizatullin $\C^*$-surfaces}
These graphs were used by Gizatullin \cite{Gi}, and systematically
studied by Dubouloz \cite{Du2}. Here we express the extended graph
of a hyperbolic Gizatullin surface $V=\Spec \C[t][D_+,D_-]$ in
terms of the divisors $D_\pm$ on $\A^1$. In \ref{complet.15} and
\ref{DGsur} we apply  these descriptions to study
Danilov-Gizatullin $\C^*$-surfaces.

\subsection{Extended graphs}
\bdefi\label{eboundary.7} Let $V$ be a Gizatullin surface and
$(\bV, D)$ a completion of $V$ by a zigzag. By Proposition
\ref{equivariant.5}(b) we can transform $(\bV, D)$ into a standard
completion $(\bV_\st, D_\st)$ so that
$$
D_\st=C_0+\ldots+C_n
$$
as in (\ref{zigzag}). We also consider the minimal resolutions of
singularities $V'$, $(\tV, D)$ and $(\tV_\st,D_\st)$ of $V$,
$(\bV, D)$ and $(\bV_\st, D_\st)$, respectively.

As in \ref{unisit} the
linear systems $|C_0|$ and $|C_1|$ define a morphism
$\Phi=\Phi_0\times\Phi_1:\tV_{\rm st}\to \PP^1\times\PP^1$ with
$\Phi_i=\Phi_{|C_i|}$, $i=0,1$. As before we choose
the coordinates in such a way that
$$C_0=\Phi_0^{-1}(\infty)\,, \qquad \Phi(C_1)= \PP^1\times
\{\infty\} \quad\mbox{and}\quad C_2\cup\ldots \cup C_n\subseteq
\Phi_0^{-1}(0) \,.$$
We recall that the divisor
$$ D_{\rm ext}=C_0\cup C_1\cup \Phi_0^{-1}(0) $$
is the {\em extended divisor} and its dual graph, also
denoted by $D_{\rm ext}$, the {\em extended graph} of
$(\tV_\st,V)$ or of $V$, for short.
\edefi

\brems\label{eboundary.8} 1. By Corollary 3.5 in \cite{FKZ} the
standard zigzag $D_{\rm st}\subseteq D_{\rm ext}$ as above is
uniquely determined up to reversing the chain $C_2,\ldots, C_n$
in
(\ref{zigzag}). However, the extended divisor $D_{\rm ext}$
usually depends on the completion.

2. As follows from Definition \ref{eboundary.7}, the extended
graph $D_{\rm ext}$ can be blown down to \be\label{3curves}
\cou{C_0}{0}\lin\cou{C_1}{0}\lin\cou{C_2}{0}\quad. \ee In
particular, $D_{\rm ext}$ is a tree, and the intersection form
$I(D_{\rm ext})$ has exactly one positive and one zero eigenvalues
(see \cite[4.1]{FKZ}).
\erems

We let $\kappa (C)$ denote the number of irreducible
components of
a curve $C$ and $\rho (V)={\rm rk} (\Pic (V))$ denote the Picard
number of $V$.

\bcor\label{picard} With $E$ being the  exceptional locus of the
minimal resolution of singularities $V'\to V$ we have
\be\label{picnum} \rho (V)=\kappa (D_{\rm ext})-\kappa (D_{\rm
st}) -\kappa (E)-1\,.\ee \ecor

\bproof Indeed,  $\Pic
(\PP^1\times\PP^1)\cong\Z[C_1]\oplus\Z[C_2]$, hence $\Pic
(\tV_{\rm st})\cong \Z^{\kappa (D_{\rm ext})-1}$ is freely
generated by the components of $D_{\rm ext}\ominus C_0$. Now
$D_{\rm st}=\tV_{\rm st}\setminus V'$, so
$$\rho (V)=\rho (V')-\rho (E)=
\rho (\tV_{\rm st})-\rho (D_{\rm st})-\rho (E) =\kappa (D_{\rm
ext})-1-\kappa (D_{\rm st})-\kappa (E)\,,$$ and the result
follows. \eproof

\bexas\label{rank} 1. For an affine  toric surface $V$ the Picard
number $\rho (V)$ vanishes, hence by (\ref{picnum}), $D_{\rm
ext}\ominus D_{\rm st}\ominus E$ consists of just one
component.

2. For $V=\PP^2\setminus C$, where $C$ is a smooth conic,
$\Pic
(V)=\Z/4\Z$ and so $\rho (V)=0$. Hence $D_{\rm ext}
\ominus D_{\rm
st}$ consists of one component, where $D_{\rm
st}=[[0,0,-2,-2,-2]]$.

3. Since any Danilov-Gizatullin surface $V_{k+1}$ (see
\ref{complet.15} below) is smooth and $\rho (V_{k+1})=1$ ,
$D_{\rm
ext}\ominus D_{\rm st}$ consists of $2$ components.

4. As was shown in Lemma\ref{pextended.01}, $V$ is toric if and
only if its extended divisor $D_{\rm ext}$ has a linear dual
graph. \eexas

\subsection{Extended graphs on Gizatullin $\C^*$-surfaces}
\bdefi\label{eboundary.4} A {\em feather} $\fF$ is a
linear chain
of smooth rational curves with dual graph \be\label{picfeather}
\co{B}\llin\boxo{e/m} \ee where $B$ has self-intersection
$\le -1$
and $e,m$ are coprime integers with $m>0$, cf. (\ref{reversed
chain}). Note that the box can also be empty. The curve $B$
will
be called the {\em bridge curve}. A collection of feathers
$\{\fF_\rho\}$ consists of feathers $\fF_\rho$, $\rho=1,\ldots
,r$, that are pairwise disjoint. Such a collection will be
denoted by a plus box $\rbh{ \xbox}\;$.

We say that this collection of feathers $\{\fF_\rho\}$ is attached
to the curve $C_i$ in a chain (\ref{zigzag}) if the bridge curves
$B_\rho$ meet $C_i$ in pairwise distinct points, and all the
feathers are disjoint with the curves $C_j$ for $j\ne i$. In a
diagram we will write in brief $$
\co{C_i}\lin\xboxo{\{\fF_\rho\}} \quad.$$ \edefi

\bexas\label{eboundary.5} 1. A curve $B$ with self-intersection
$\le -1$ is a feather; the box in (\ref{picfeather}) is in this
case empty and there is only a bridge curve.

2. The {\em $A_k$-feather} is
$$
\cou{E}{-1}\lin \boxo{A_k}\quad.
$$
Thus the $A_k$-feather represents the contractible linear chain
$[[-1,(-2)_{k}]]$ and the $A_0$-feather represents a single
$(-1)$-curve $E$. \eexas

\bdefi\label{eboundary.6} A feather collection will be called
{\em
admissible} if it contains at most one feather which is not an
$A_k$-feather. \edefi

In the next proposition we describe the extended graphs of
Gizatullin surfaces with a hyperbolic $\C^*$-action. We recall
that for such a surface $V=\Spec A_0[D_+,D_-]$, necessarily
$A_0=\C[t]$ and both $\supp (\{D_\pm\})$ consist of at most one
point.

\bprop\label{eboundary.9}
\begin{enumerate} [(a)]
\item The minimal resolution $V'\to V$ of a non-toric Gizatullin
surface $V$ with a hyperbolic $\C^*$-action admits a standard
equivariant completion $(\tV_\st, D_\st)$ with $D_{\rm
st}=C_0+\ldots +C_n$ as in Definition \ref{eboundary.8} and with
the following extended graph $D_{\rm ext}$ :

\bigskip

\be\label{ezigzag1} \cou{C_0}{0}\lin\cou{C_1}{0}\lin
\cou{C_2}{w_2}\lin\ldots
\lin\cou{}{w_s}\nlin\xbshiftup{}{\qquad \{\fF_\rho\}_{\rho \ge
1}}\lin \cou{C_{s+1}}{w_{s+1}}\lin \ldots \lin \cou{C_n}{w_n}
\llin\xboxo{\fF_0} \quad,\ee where $w_i\le -2$ $\forall i\ge 2$,
$\fF_0$ is a single feather and $\{\fF_\rho\}_{\rho\ge 1}$ is a
nonempty admissible feather collection. \item If, moreover, $\supp
(\{D_+\})\cup \supp (\{D_-\})$ consists of at most one point then,
after possibly reversing the chain $(C_2,\ldots, C_n)$ in the
standard zigzag, we can achieve additionally that
\begin{enumerate}[(i)]
\item the chain \qquad
$\co{C_{s+1}}\lin\ldots\lin\co{C_n}\lin\xboxo{\fF_0}$ \qquad is
either empty or is not contractible to a smooth point, and

\item all the $\fF_\rho$, $\rho\ge 1$, are $A_{s_\rho}$-feathers
for some $s_\rho\ge 0$.
\end{enumerate}
\item If $\supp (\{D_+\})\cup \supp (\{D_-\})$ consists of two
points then the chain in (i) contracts to a smooth point.
\end{enumerate}\eprop

\bproof Let as before $V=\Spec A_0[D_+,D_-]$ be a DPD-presentation
of $V$ with $A_0=\C[t]$. Since $V$ is a Gizatullin surface, we
have $\supp (\{D_\pm\}) \subseteq \{p_\pm\}$ for some points
$p_+,p_-\in \PP^1$. So $p_+,p_-$ are among the points
$\{p_i,q_j\}$ considered in \ref{complet.7}, and $\supp
(\{D_\pm\})$ can also be empty or equal. We will construct a
standard equivariant completion $(\tV_\st, D_\st)$ of $V'$
starting from the natural completion $(\tV,\tD)$ as obtained in
\ref{complet.12}.

\noindent $\bullet$ Let us first consider the case where
$p_+=p_-$. In this case, after passing to an equivalent pair
$(D_+,D_-)$ if necessary, none of the $q_j$ is present besides
possibly $p_+$, and for all the $p_i$ different from $p_+$ the
numbers $D_\pm(p_i)$ are integral. According to Example
\ref{complet.13}, the fiber of $\tilde\pi:\tV\to \PP^1$ over
$p_i\ne p_+$ together with the sections $\tC_\pm$ of $\pi$ is
$$
\co{\tC_+}\lin \co{\tO_i^+}\lin\boxo{A_{s_i}}\lin \co{\tO_i^-}
\llin\co{\tC_-} $$ with $s_i=-1-(D_+(p_i)+D_-(p_i))$. The fiber
$\tilde\pi^{-1}(p_+)$ together with the sections $\tC_\pm$ is in
case $D_+(p_+)+D_-(p_+)= 0$ $$
\co{\tC_+}\vlin{15}\boxo{\{D_+(p_+)\}}\vlin{15}\co
{\tO_{p_+}}\vlin{15}\boxo{\{D_-(p_+)\}^*}\vlin{15}\co{\tC_-}
\quad,$$ and in case $D_+(p_+)+D_-(p_+)< 0$ \be\label{afiber}
\co{\tC_+}\vlin{15}\boxo{\{D_+(p_+)\}}\vlin{15} \co
{\tO_{p_+}^+}\lin\boxo{R_{p_+}}\lin\co{\tO_{p_+}^-}
\vlin{15}\boxo{\{D_-(p_+)\}^*}\vlin{15}\co{\tC_-}\quad, \ee where
$R_{p_+}$ stands for the minimal resolution of the cyclic quotient
singularity in the fiber $\pi^{-1}(p_+)$, see Proposition
\ref{complet.12}(c). By Corollary \ref{boundarydivisor}, in both
cases the boundary zigzag is
  \be\label{eqbdiv}\quad
\boxo{\{D_+(p_+)\}^*}\vlin{13}\co{\tC_+}\lin \co{F_\infty}\lin
\co{\tC_-} \vlin{13}\boxo{\{D_-(p_-)\}}\quad\quad,\ee where
$p_+=p_-$, $F_\infty^2=0$ and, according to Proposition
\ref{complet.12}(d), \be\label{vertexweight} \tC_+^2+\tC_-^2=\deg
\left(\lfloor D_+\rfloor + \lfloor D_-\rfloor\right)\le \deg
\left( D_++ D_-\right)\le 0\,. \ee

\smallskip

{\it Claim. (a) If $\tC_+^2+\tC_-^2\ge -1$ then $V$ is a toric
surface. (b) If $\tC_+^2+\tC_-^2\le -2$ then moving the zero
weight of $F_{\infty}$ in (\ref{eqbdiv}) to the left yields a
standard zigzag.}

\smallskip

{\it Proof of the claim.} (a) If $\tC_+^2+\tC_-^2=0$ then by
virtue of (\ref{vertexweight}), $D_+=-D_-$ and both divisors are
integral. Hence the pair $(D_+,D_-)$ is equivalent to $(0,0)$ and
$V\cong \A^2$ is toric. If $\tC_+^2+\tC_-^2=-1$ then either
$D_+=-D_-$ and $D_\pm(a)$ are integers except for one point $a=p$,
or there is a point $p\in\A^1$ where $D_+(p)+ D_-(p)<0$, and for
all the other points $q\in\A^1$, $q\neq p$, we have $D_+(q)+
D_-(q)=0$ and $\{D_\pm (q)\}=0$. Anyhow, passing to an equivalent
pair of divisors $(D_+,D_-)$ we obtain
$$D_+=-\frac{e_+}{m_+}[p]\qquad\mbox{and}\qquad D_-=
\frac{e_-}{m_-}[p]\,.$$
By Lemma \ref{torics}(b), in this case $V$ is toric. This shows
(a).

To show (b) we perform inner (hence equivariant) elementary
transformations in (\ref{eqbdiv}) which replace the curves $\tC_+$
and $F_{\infty}$ by two others with self-intersection $0$ making
the new weight of $\tC_-$ equal to $\tC_+^2+\tC_-^2\le -2$.
Further we perform inner elementary transformations moving the
two zeros to the left to obtain the boundary zigzag on $\tV_{\rm
st}$ in the standard form
$$\quad \cou{C_0}{0}\lin
\cou{C_1}{0}\vlin{13}\boxo{\{D_+(p_+)\}^*}\vlin{13}
\cou{\tC_-}{\le-2} \vlin{13}\boxo{\{D_-(p_+)\}}\quad,$$ see Lemma
2.12 in \cite{FKZ}. These elementary transformations do not
contract the components to the right of $\tC_-$ preserving their
weights. \qed

\smallskip

Since by our assumption the surface $V$ is non-toric, we are in
case (b) above. We attach to the curve $\tC_-$ the collection of
feathers \quad$\fF_i:\,\boxo{A_{s_i}}\lin
\co{\tO_i^-}\,\,\,,\,\,\,$ and in case $D_+(p_+)+D_-(p_+)\ne 0$ to
the last curve of the weighted $\{D_-(p_+)\}$-box also the feather
\be\label{0feather} \fF_0:\qquad
\boxo{R_{p_+}}\lin\co{\tO_{p_+}^-}\quad. \ee This leads to the
graph
\bigskip

\be\label{eqegraph} \cou{C_0}{0}\lin
\cou{C_1}{0}\vlin{15}\boxo{\{D_+(p_+)\}^*}\vlin{13}
\cou{}{\hC_-}\nlin\xbshiftup{}{\qquad \{\fF_\rho\}_{\rho \ge 1}}
\vlin{15} \boxo{\{D_-(p_+)\}} \vlin{15}\xboxo{\fF_0} \quad,\ee
where $\hC_-$ is the proper transform of $\tC_-$ with $\hC_-^2\le
-2$. We claim that (\ref{eqegraph}) is already the full extended
graph $D_{\rm ext}$ or, equivalently, that the curves in
(\ref{eqegraph}) besides $C_0,C_1$ constitute the full fiber
$\Phi_0^{-1}(0)$.

In fact, all the components of $D_{\rm ext}$ are $\C^*$-stable,
since so are the curves $C_0,C_1$ and the linear systems $|C_0|$
and $|C_1|$ on $\tV_{\rm st}$. Moreover, since the extended graph
is a tree, a curve which occurs in $\Phi_0^{-1}(0)\ominus D_{\rm
st}$ meets the boundary zigzag $D_{\rm st}$ in at most one point.
Thus the proper transforms on $\tV_{\rm st}$ of the curves
$\tO^+_{p_i}$ and $\tO^+_{p_+}$, respectively, $\tO_{p_+}$ or of
an irreducible fiber of $\tilde \pi:\tV\to\PP^1$ cannot appear in
$\Phi_0^{-1}(0)$. All the other $\C^*$-invariant curves belong
already to the boundary zigzag $D_{\rm st}$ or are in one of the
feathers (indeed, in $\tV$ the only $\C^*$-invariant curves are
those in the fibers and the curves $\tC_\pm$), proving the claim.

Now (ii) is clear from the construction. To deduce (i), assume in
contrary that the chain in (i) is contractible to a smooth point.
This is only possible in the case where $D_+(p_+)+D_-(p_+)< 0$,
since otherwise the feather $\fF_0$ is empty by construction.
Moreover, $\tO_{p_+}^-$ in the feather $\fF_0$ in (\ref{0feather})
must be a $(-1)$-curve, since otherwise the chain in (i) would be
minimal, contrary to our assumption. Thus as well the part
$$
P_{p_+}:\qquad \boxo{R_{p_+}}\lin\co{\tO_{p_+}^-}\vlin{17}
\boxo{\{D_-(p_+)\}^*}
$$
of the fiber $\tilde\pi^{-1}(p_+)$ in (\ref{afiber}) is
contractible to a smooth point. After contracting this, the
resulting dual graph of the blown down fiber must be \quad\qquad
$\boxo{\{D_+(p_+)\}^*}\vlin{15}\co{-1}\,\,\,$, or it is just
irreducible. The former case is not possible since obviously a
linear graph with one $(-1)$-end vertex and the rest of vertices of
weights at most $-2$, cannot be contracted to a single vertex of
weight $0$.
In the latter case $D_+(p_+)$ is integral i.e., $m_+=1$. If we now
interchange $D_+$ and $D_-$ then the chain in (i) is either not
contractible, or, by the same argument, in addition $D_-(p_+)=1$.
In the latter case the chain in (i) is empty, as desired, see
(\ref{eqegraph}). In case $p_+=p_-$ this proves (a) and (b).

\noindent $\bullet$ Suppose further that $p_+$ and $p_-$ are
distinct, so that $\{D_-(p_+)\}=\{D_+(p_-)\}=0$. With the same
reasoning as before, the completion $\tV$ from Proposition
\ref{complet.11} has the boundary divisor as in (\ref{eqbdiv}),
and $\tV_{\rm st}$ has the extended graph as in (\ref{eqegraph})
with the only difference that now one of the
$\{\fF_\rho\}_{\rho\ge 1}$ might be not an $A_s$-feather. We leave
the details to the reader. This completes the proof of (a).

To show (c) we note that by our construction, up to reversing,
the chain in (i) is the part
$$
P_{p_-}:\qquad \boxo{R_{p_-}}\lin\co{\tO_{p_-}^-}\vlin{17}
\boxo{\{D_-(p_-)\}^*}
$$
of the linear chain $\tilde\pi^{-1}(p_-)+\tC_++\tC_-$ :
$$
\co{\tC_+}\llin\co{\tO_{p_-}^+}\lin\boxo{R_{p_-}}
\lin\co{\tO_{p_-}^-}
\vlin{15}\boxo{\{D_-(p_-)\}^*}\vlin{15}\co{\tC_-}\quad,
$$
where $\{D_-(p_-)\}^*$ is associated to the cyclic quotient
singularity of type $(-m^-,-e^-)$. We claim that $\tO_{p_-}^-$ is
a $(-1)$-curve. Indeed, if $\tO_{p_-}^-$ were not a $(-1)$-curve
then $\tO_{p_-}^+$ would be since at least one of these is a
$(-1)$-curve by Proposition \ref{complet.9}. But then the fiber
$\tilde\pi^{-1}(p_-)$ would represent a linear chain with just one
$(-1)$-curve at the end and so it cannot be contracted to a
$0$-curve by the same argument as before. Now tracing the process
of blowing down to a non-degenerate fiber $[[0]]$, it is easily
seen that $P_{p_-}$ is indeed contractible to a smooth point.
\eproof

\brem\label{eboundary.10} The curve on the left of $\hC_-$ in
(\ref{eqegraph}) has self-intersection $-2$ if and only if
$m^+/e_+\le 2$, where $e_+$ is the integer with $0\le e_+< m^+$
and $e_+\equiv -e^+\mod m^+$. Similarly, the curve on the right
of $\hC_-$ in (\ref{eqegraph}) has self-intersection $-2$ if and
only if $-m^-/e_-\le 2$, where $e_-$ is the unique integer with
$0\le e_-< -m^-$ and $e_-\equiv -e^-\mod m^-$.\erem

The following result is a useful supplement to Corollary
\ref{corzigzag}.

\bcor\label{ab} If the surface $V$ in Proposition
\ref{eboundary.9} is smooth then the following hold.
\begin{enumerate} [(a)] \item In case where
$\supp (\{D_+\})\cup\supp (\{D_-\})=\{p_0\}$, up to reversing the
zigzag, the extended graph $D_{\rm ext}$ of $V$ is

\bigskip

$$
\bigskip
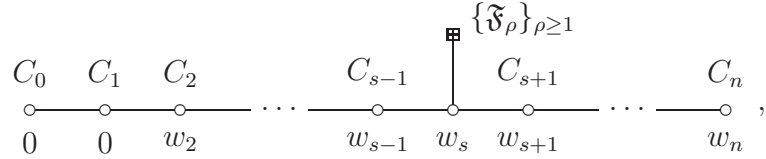
\cou{C_0}{0}\lin
\cou{C_1}{0}\lin\cou{C_2}{w_2}\lin\ldots\lin\cou{C_{s-1}}
{w_{s-1}}\lin
\cou{}{w_s}\nlin\xbshiftup{}{\qquad \{\fF_\rho\}_{\rho \ge 1}}
\lin \cou{C_{s+1}}{w_{s+1}}\lin\ldots \lin\cou{C_{n}}{w_n}
\quad,$$ where $w_i\le -2$ $\forall i\ge 2$, every feather
$\fF_\rho$ consists of a single $(-1)$-vertex and the number of
these feathers is equal to $|w_s|-1$. Moreover, the length
$\kappa
(D_{\rm st})$ of the boundary zigzag $D_{\rm st}$ can be $n+1$ or
$n$ depending on whether $C_n$ is in $D_{\rm st}$ or not, and
$$\rho(V)={\rm rk} (\Pic (V))=\begin{cases}
|w_s|-2 & \,\,\,\mbox{if}\,\,\,\kappa (D_{\rm st})=n+1\\
|w_s|-1 & \,\,\,\mbox{if}\,\,\,\kappa (D_{\rm
st})=n\,.\end{cases}$$ \item In case where $\supp
(\{D_\pm\})=\{p_\pm\}$ with $p_+\ne p_-$, up to
equivalence of the
pair $(D_+,D_-)$ we have $D_+(p_+)=-1/m_+$, $D_-(p_+)=0$ and
$D_+(p_-)=0$, $D_-(p_-)=-1/m_-$ with $m_+,m_-\ge 2$.
The extended
graph $D_{\rm ext}$ of $V$ is

\bigskip

$$
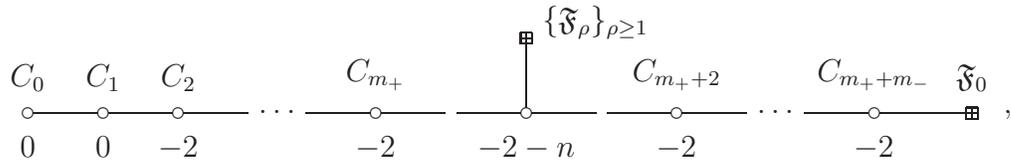
\cou{C_0}{0}\lin
\cou{C_1}{0}\lin\cou{C_2}{-2}\lin\ldots\lin\cou{C_{m_+}}
{-2}\lin\lin
\cou{}{-2-n}\nlin\xbshiftup{}{\qquad \{\fF_\rho\}_{\rho \ge 1}}
\lin\lin \cou{C_{m_++2}}{-2}\lin\ldots \lin\cou{C_{m_++m_-}}{-2}
\llin\xboxo{\fF_0} \quad,$$ where $\fF_1$ is a feather consisting
of a single $(-m_+)$-curve $\tO_{p_+}^-$, $\fF_\rho$, $\rho> 1$,
are $n$ feathers consisting of $(-1)$-curves $\tO_{p_\rho}^-$ and
$\fF_0$ is a feather consisting of a single $(-1)$-curve
$\tO_{p_-}^-$.
\end{enumerate}\ecor

\bexas\label{6.9} 1. It is clear that the curves $\tC_\pm$ in the
completion constructed in Proposition \ref{complet.12} are
pointwise fixed by the $\C^*$-action. Thus the component $C_s$ as
in Proposition \ref{eboundary.9} joined by bridges with a feather
collection is parabolic.

2.  In the case when $D_+=0$ and $D_-=-2/3[a]$ for some
$a\in\A^1$, we let $\tV$ be the resolution of the completion of
$V=\Spec A_0[D_+,D_-]$ constructed in Proposition
\ref{complet.6}. Then $V$ is a toric surface and its boundary in
$\tV$ has dual graph
$$
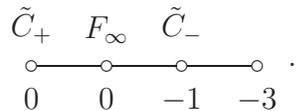
\cou{\tC_+}{0}\lin \cou{F_\infty}{0}\lin \cou{\tC_-}{-1}\lin
\cou{}{-3}\quad.
$$
Blowing up the intersection point $\tC_+\cap F_\infty$ and
contracting the proper transforms of $\tC_\pm$ leads to an
equivariant  completion of $V$ by a standard zigzag $[[0,0,-2]]$
and without any $\C^*$-parabolic component. \eexas

The latter cannot happen for a non-toric $\C^*$-surface, see Lemma
\ref{parabolic}.

In the following result we analyze to what extent the
extended graph determines a non-toric Gizatullin $\C^*$-surface.

\bprop\label{eboundary.thm} Suppose that two non-toric Gizatullin
$\C^*$-surfaces have the same extended graphs and the same
positions of the feathers on the parabolic component. Then these
surfaces are equivariantly isomorphic.
\eprop

\bproof This can be easily derived from the fact that the
DPD-presentation determines the $\C^*$-surface uniquely up to an
equivariant isomorphism. \eproof

\subsection{Danilov-Gizatullin $\C^*$-surfaces}
The following class of examples was elaborated by Danilov and
Gizatullin \cite{DaGi} (see also the Introduction). Answering our
question on the uniqueness of $\C^*$-actions \cite{FlZa2}, P.\
Russell showed that there are several non-conjugated
$\C^*$-actions on a Danilov-Gizatullin surface. We expose here
these $\C^*$-actions in a somewhat different manner.

\bexa\label{complet.15} Given a pair of natural numbers $k,r$ with
$1\le r\le k$ and a pair of distinct points
$p_0,\,p_1\in\A^1=\Spec \C[t]$, we consider the smooth affine
hyperbolic $\C^*$-surface $V=V_{k,r}=\Spec A_0[D_+,D_-]$, where
$A_0=\C[t]$, \be\label{GDS}
D_+=-\frac{1}{r}[p_0]\quad\mbox{and}\quad
D_-=-\frac{1}{k+1-r}[p_1]\,.\ee We call these {\it
Danilov-Gizatullin $\C^*$-surfaces}.

By Lemma \ref{singtype}, the equivariant completion $\bV$ of $V$
as constructed in Proposition \ref{complet.6} has an
$A_{r-1}$-singularity at the point $p_0^+$ and an
$A_{k-r}$-singularity at $p_1^-$, whereas the other points shown
at the following diagram are smooth:

\bigskip

$$
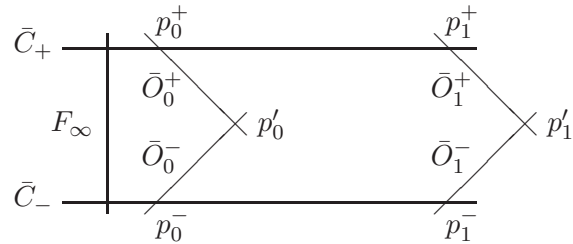
\unitlength0.5cm%
\begin{picture}(12,5.0)
\small \put(1.3,0.8){\line(1,0){11}} \put(1.3,4.9){\line(1,0){11}}
 \put(2.5,0.5){\line(0,1){4.8}}

\put(3.5,5.3){\line(1,-1){2.7}} \put(3.5,0.5){\line(1,1){2.7}}
\put(11.2,5.3){\line(1,-1){2.7}} \put(11.2,0.5){\line(1,1){2.7}}

\put(3.8,5.4){$p_0^+$} \put(3.8,0){$p_0^-$}
\put(6.5,2.7){$p_0'$} \put(3.4,3.6){$\bO_0^+$}
\put(3.4,1.8){$\bO_0^-$} \put(1.0,2.7){$F_\infty$}

\put(11.5,5.4){$p_1^+$} \put(11.5,0){$p_1^-$}
\put(14.2,2.7){$p_1'$} \put(11.1,3.6){$\bO_1^+$}
\put(11.1,1.8){$\bO_1^-$}

\put(0.0,4.8){$\bC_+$} \put(0.0,0.7){$\bC_-$}
\end{picture}
$$

\medskip

\noindent here  $\bO_0^\pm:=\bO_{p_0}^\pm$ and
$\bO_1^\pm:=\bO_{p_1}^\pm$. By Corollary \ref{boundarydivisor},
the
boundary zigzag $\bar D\subseteq\bar V$ is
$$
\boxo{A_{r-1}}\lin\cou{\tC_+}{-1}\lin \cou{F_\infty}{0}
\lin\cou{\tC_-}{-1}\lin \boxo{A_{k-r}}\quad,
$$ where $F_\infty$ denotes the fiber of $\pi$ over $\infty\in
\PP^1$. Contracting successively all $(-1)$-curves provides an
equivariant completion $\bV_{k,r}$ of $V_{k,r}:=V$ by a single
smooth rational curve, say, $S$ of self-intersection $k+1$. For a
fixed $k$, by a theorem of Danilov-Gizatullin \cite{DaGi} the $k$
affine surfaces $V_{k,r}$, $1\le r\le k$, are all isomorphic.
However, by Theorem 4.3(b) in \cite[I]{FlZa1} they are not
equivariantly isomorphic since the fractional parts of the pairs
$(D_+,D_-)$ are all distinct for distinct values of $r$. Thus the
Danilov-Gizatullin surface $V_{k+1}\cong V_{k,r}$ possesses at
least $k$ different $\C^*$-actions that are not conjugated in the
automorphism group $\Aut (V_{k+1})$. Furthermore the action of the
automorphism $\lambda\longmapsto\lambda^{-1}$ of the group $\C^*$
amounts to interchanging $D_+$ and $D_-$, which reduces the number
of essentially different $\C^*$-structures on $V_{k+1}$ to
$\lfloor\frac{k+1}{2}\rfloor$.

Let us study the extended divisors $D_{\rm ext}$ of the
$\C^*$-surfaces $V_{k,r}\cong V_{k+1}$. The fibers over $p_0$ and
$p_1$ together with the curves $\tC_\pm$ have dual graphs
$$
\cou{\tC_+}{-1}\lin\boxo{A_{r-1}}\lin \cou{\bO_0^+}{-1} \lin
\cou{\bO_0^-}{-r} \lin\cou{\tC_-}{-1}\,\,\,,
\qquad\mbox{respectively,}\qquad \cou{\tC_+}{-1}\vlin{15}
\cou{\bO_1^+}{-(k+1-r)} \vlin{15} \cou{\bO_1^-}{-1}
\lin\boxo{A_{k-r}} \lin\cou{\tC_-}{-1}\quad .
$$
Thus moving the zero on the
boundary to the left by means of elementary transformations leads
to the extended graph \vspace*{1truecm}
$$
\cou{C_0}{0}\lin\cou{C_1}{0}\lin
\boxo{A_{r-1}}\lin\cou{\qquad\bC_-}{-2} \nlin\cshiftup{\bO_0^-}{-r}
\llin \boxo{A_{k-r}}\lin\cou{\bO_1^-}{-1}\quad ,
$$
where the feathers are formed by $\bO_0^-$ and $\bO_1^-$.
Similarly,
moving the zero to the right leads to the extended graph
\vspace*{1truecm}
$$
\cou{\bO_0^+}{-1}\lin \boxo{A_{r-1}}\lin\cou{\qquad\bC_+}{-2}
\nlin\cshiftup{\bO_1^+}{\qquad \qquad -(k-r+1)} \llin
\boxo{A_{k-r}}\lin\cou{C_1}{0}\lin\cou{C_0}{0}\quad ,
$$
where now the feathers are $\bO_0^+$ and $\bO_1^+$. In both cases,
the standard
boundary zigzag $D_{\rm st}$ is $[[0,0,(-2)_k]]$ with dual graph
$$\cou{C_0}{0}\lin\cou{C_1}{0}\lin \boxo{A_{k}}\quad.$$
\eexa

\bprop\label{DGsur} The Danilov-Gizatullin surface $V_{k+1}$
($k\ge 0$) carries exactly $k$ different, up to conjugation
in the
automorphism group, $\C^*$-actions, and all of them are
hyperbolic. \eprop

Let us give two alternative proofs.

{\it 1-st proof.} A smooth elliptic or parabolic Gizatullin
$\C^*$-surface is necessarily isomorphic to $\A^2$, see Corollary
4.4 in \cite[II]{FlZa1}. Hence the Gizatullin surface $V_{k+1}$
with the Picard group $\Pic (V_{k+1})\cong \Z$ cannot carry any
elliptic or parabolic $\C^*$-action.

We have shown in \ref{complet.15} above that there are at least
$k$ mutually non-conjugated hyperbolic $\C^*$-actions on
$V_{k+1}$. To show that any such action on $V_{k+1}$ is conjugated
to one of these is the same as to show that, given an isomorphism
$$V_{k+1}\cong\Spec A_0[D_+,D_-]$$ with $A_0=\C[t]$ and some pair
of $\Q$-divisors $D_\pm$ on $\A^1=\Spec
 A_0$ with
$D_++D_-\le 0$, up to equivalence $(D_+,D_-)$ must be one of the
pairs (\ref{GDS}). Since $V_{k+1}$ is a Gizatullin surface, the
supports of $\{D_+\}$ and $\{D_-\}$ consist of at most one point.
Let as before $p_0,\ldots,p_l$ be the points with
$D_+(p_i)+D_-(p_i)<0$, and $q_1,\ldots,q_s$ the points with
$D_+(q_j)+D_-(q_j)=0$. Replacing $D_+,\,D_-$ by an equivalent pair
we may suppose that $\{D_\pm(q_j)\}\neq 0$. Thus necessarily $s\le
1$. If $s=1$ then by Corollary 4.24 in \cite[I]{FlZa1}, $\Pic
(V_{k+1})$ would have torsion. Since $\Pic (V_{k+1})\cong \Z$,
this case is impossible and so $s=0$.
On the other hand, again by Corollary 4.24 in \cite[I]{FlZa1}, we
have $l=1$.

First we assume that both $\{D_+(p_0)\}$ and
$\{D_-(p_0)\}$ are nonzero. Then necessarily
$\{D_+(p_1)\}=\{D_-(p_1)\}=0$.
As $p_0'\in V_{k+1}$ is a smooth point, by Lemma
\ref{singtype}(c) we have
$$
D_+(p_0)+D_-(p_0)=\frac{\Delta (p_0)}{m_0^+m_0^-}
=\frac{1}{m_0^+m_0^-}\,\,.
$$
This implies $\lfloor D_+(p_0)\rfloor+\lfloor D_-(p_0)\rfloor=-1$.
The standard boundary zigzag of $V_{k+1}$ is \be\label{staba}
\cou{C_0}{0}\lin
\cou{C_1}{0}\vlin{15}\boxo{\{D_+(p_0)\}^*}\vlin{13} \cou{}{\omega}
\vlin{15} \boxo{\{D_-(p_0)\}} \qquad\quad = \quad
[[0,0,(-2)_k]]\,,\ee where $\omega=\sum_{i=0,1} \left(\lfloor
D_+(p_i) \rfloor+\lfloor D_-(p_i) \rfloor\right)=D_+(p_1) +
D_-(p_1) -1=-2$. Therefore $D_+(p_1)+ D_-(p_1)=-1$. Moreover, the
boxes in (\ref{staba}) are $A_{r-1}$ and $A_{k-r}$-boxes for some
$r$ with $0<r<k+1$, so that \be\label{0frac}
\{D_+(p_0)\}=\frac{r-1}{r}\qquad\mbox{and}\qquad
\{D_-(p_1)\}=\frac{k-r}{k+1-r}\,.\quad \ee Passing to an
equivalent pair of divisors we may assume that
$D_+(p_0)=\frac{-1}{r}$, hence $\lfloor D_+(p_0)\rfloor = -1$,
$\lfloor D_-(p_0)\rfloor = 0$ and $D_-(p_0) =
\frac{k-r}{k+1-r}=\frac{e_0^-}{m_0^-}$, where $e_0^-=-(k-r)$ and
$m_0^-=-(k+1-r)$. Again by smoothness of the point $p_0'\in
V_{k+1}$, the determinant (\ref{dete}) is equal to $1$:
$$\Delta (p_0)=-\left|\ba{ll} 1 & -(k-r)\\ r &
-(k+1-r)\ea\right|=1\,.$$ Hence $(k+1-r)-(k-r)r=1$ and so,
$(k-r)(r-1)=0$. This forces $k=r$ or $r=1$ i.e., $D_+(p_0)$ or
$D_-(p_0)$ is integral, contrary to our assumption.

Therefore, up to interchanging $p_0$ and $p_1$, the only
possibility is
  $$\{D_+(p_0)\}\neq 0\qquad\mbox{and}\qquad \{D_-(p_1)\}\neq
0\,,$$ whereas $D_-(p_0)$ and $D_+(p_1)$ are integral. After
passing again to an equivalent pair $(D_+,D_-)$ we may suppose
that $D_-(p_0)=D_+(p_1)=0$. We write now
$$D_+ (p_0)=-\frac{e_0}{m_0}\qquad\mbox{and}\qquad
D_- (p_1)=-\frac{e_1}{m_1}\,,$$ where $m_0,\,m_1>0$.
By smoothness
of the points $p_i'\in V_{k+1}$  we have $\Delta (p_i)=1$
for $i=0,1$,
hence $e_0=e_1=1$. Thus the zigzag (\ref{staba}) of the
equivariant standard completion $(\tV_{k+1})_{\rm st}$ of
$V_{k+1}$ is
$$\cou{C_0}{0}\lin
\cou{C_1}{0}\vlin{15}\boxo{A_{m_0-1}}\vlin{13} \cou{}{-2}
\vlin{15} \boxo{A_{m_1-1}} \quad.$$ This yields $ m_0+m_1=k-1$,
so
$(D_+,D_-)$ is one of the pairs in (\ref{GDS}), as required.
  \qed

{\it 2-nd proof.} We must show that any hyperbolic $\C^*$-action
$\Lambda$ on $V_{k+1}$ is conjugate to one of those constructed in
Example \ref{complet.15}. Since these $k$ $\C^*$-actions on $V$
are mutually non-conjugate, this would complete the proof.

The surface $V_{k+1}$ admits an equivariant completion
$(\bV_{k+1})_{\rm st}$ by a standard zigzag $D=C_0+C_1+ \ldots +
C_{k+1}$ such that $C_j^2=-2$ for $j\geq 2$. As before, the
complete linear systems $|C_0|$ and $|C_1|$ yield a morphism $\Phi
=(\Phi_0,\Phi_1) : (\bV_{k+1})_{\rm st} \to Q=\PP^1 \times \PP^1$.
Since the $\C^*$-action on $(\bV_{k+1})_{\rm st}$ stabilizes $C_0$
and $C_1$ it preserves the corresponding linear systems and hence
induces a linear $\C^*$-action $(x,y) \to (\lambda^nx,
\lambda^my)$ on $Q$ such that $\Phi$ is equivariant. We note that
the numbers $n$ and $m$ uniquely determine the part of the
extended graph $D_{\rm ext}$ between $C_2$ and the parabolic
component $C_r$. Indeed $C_r$ appears as the $(-1)$-curve in the
resolution graph $\G_0$ of the curve singularity $x^m=y^n$. Unless
$n=r-2$ and $m=r-1$ for some $r$ the part of $\G_0$ between $C_2$
and $C_r$ contains vertices of weight $\leq -3$ which contradicts
our assumption. Thus $n=r-2$ and $m=r-1$ and so, besides $C_3,
\ldots ,C_r$, $\G_0$ must contain an extra vertex $E$ of weight
$E^2=-r$ which is the proper transform of a feather (unless,
maybe, in the case where $r=3$). The only way to get $C_j^2=-2$
for $j \geq r$ is to construct a linear chain $C_{r+1}, \ldots ,
C_{k+1}, E'$ with $C_j^2=-2$, where $E'$ with $(E')^2 =-1$ is the
second feather and a neighbor of $C_{k+1}$. This produces exactly
the same extended graph

\bigskip

$$
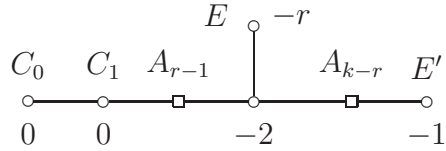
\cou{C_0}{0}\lin\cou{C_1}{0}\lin \boxo{A_{r-1}}\lin\cou{}{-2}
\nlin\cshiftup{E}{-r} \llin \boxo{A_{k-r}}\lin\cou{E'}{-1}\quad
$$ \medskip

\noindent as in Example \ref{complet.15} i.e., the same extended
graph as one of the standard actions. Now the pair of divisors
$(D_+,D_-)$ can be read up from this graph and so it coincides
with the corresponding pair (\ref{GDS}). Hence the corresponding
$\C^*$-actions on $V_{k+1}$ are conjugate. \qed

\brems\label{plusac} (1) Every Gizatullin surface $V$ admits at
least two different affine rulings (that is, $\A^1$-fibrations)
$v_\pm : V \to \A^1$. They are provided by the projections
$\Phi_0^\pm: \bV_{\rm st}^\pm\to \PP^1$ as in Definition
\ref{eboundary.7}, where $(\bV_{\rm st}^\pm, D_{\rm st}^\pm)$ are
two equivariant completions of $V$ by standard zigzags $D_{\rm
st}^+,\, D_{\rm st}^-$ which differ by reversion moving a pair
of zeros from the left to the right. By Lemma \ref{unilem}
$v_\pm : V \to \A^1$ is a smooth $\A^1$-fibration over
$\A^1\setminus \{0\}$.

If moreover $V=\Spec A_o[D_+,D_-]$ is not toric and carries a
$\C^*$-action then taking the two unique equivariant completions
there are unique affine rulings $v_\pm$ that are equivariant with
respect to a suitable $\C^*$-action on $\A^1$. Moreover we can
describe their fibers over 0 in terms of $D_\pm$: they are
disjoint unions of the $\C^*$-orbit closures $\bar O_i^\mp\cong
\A^1$, one for each point $p_i\in \A^1_\C$ with
$(D_++D_-)(p_i)<0$, see Proposition 3.25 in \cite{FlZa1}, where
also the multiplicity of $\bar O_i^\mp$ in $\div \,(v_\pm)$ is
computed.

(2) In particular, to any given hyperbolic $\C^*$-action on a
Danilov-Gizatullin surface $V_{k+1}$  corresponds  such a unique
pair $v_\pm$ of equivariant affine rulings $V_{k+1}\to \A^1_\C$.
Given $r$ with $1\le r \le k$ as in Example \ref{complet.15}
above, $v_\pm^{-1}(0)$ consists of the corresponding feather
components $\bar O_0^\mp, \bar O_1^\mp$, with multiplicities
$$\div \,(v_+)=[\bar O_0^-]+r[\bar O_1^-]\quad\mbox{and}\quad
\div \,(v_-)=
(k+r-1)[\bar O_0^+]+[\bar O_1^+]\,,
$$
see Proposition 3.25 in \cite{FlZa1}. Alternatively, these
equalities can be seen following the construction of the feather
components in the second proof above. Since conjugate
affine rulings must have equal sequences of multiplicities of
degenerate fibers, we obtain the following corollary\footnote{This
was first proved by Peter Russell, see \cite{CNR}, and in some
particular cases by Adrien Dubouloz.}. \erems

\bcor\label{cplu} (a) The Danilov-Gizatullin surface $V_{k+1}$
($k\ge 0$) carries at least $\lfloor\frac{k+1}{2}\rfloor$
different, up to conjugation in the automorphism group, affine
rulings $V_{k+1}\to \A^1_\C$ with a unique degenerate
fiber.

(b) Given integer $r\neq \frac{k+1}{2}$ with $1\le r \le k$, the
equivariant affine rulings $v_\pm: V_{k,r}=V_{k+1}\to \A^1_\C$
canonically attached to the corresponding $\C^*$-action on
$V_{k+1}$ are not conjugate. \ecor

See also \cite{Du3} for another approach to (a).


\begin{thebibliography}{KMMH}

\bibitem[BML]{BML}
T.\ Bandman, L.\ Makar-Limanov, {\em Affine surfaces with
$AK(S)=\C$}. Michigan Math.\ J.\ 49 (2001), 567--582.

\bibitem[BH]{BH} F. Berchtold, J. Hausen, {\em Demushkin's
theorem in codimension one}, Math.\ Z.\ 244 (2003), 697--703.

\bibitem[Be]{Be} J.\ Bertin,
{\em  Pinceaux de droites et automorphismes des surfaces affines}.
J.\ Reine Angew.\ Math.\ 341 (1983), 32--53.

\bibitem[CNR]{CNR} P.\ Cassou-Nogu\`es, P.\ Russell,
{\em Birational endomorphisms of the affine plane and surfaces of
the form "rational ruled surface minus a section"} ({\it in
progress}).

\bibitem[DaGi]{DaGi}
V.\ I.\ Danilov, M.\ H.\ Gizatullin, {\em Automorphisms of affine
surfaces}. I. Math. USSR Izv. 9 (1975), 493--534; II. {\it ibid.}
11 (1977), 51--98.

\bibitem[DaiRu]{DaiRu}
D.\ Daigle, P.\ Russell, I. {\em Affine rulings of normal rational
surfaces}. Osaka J.\ Math.\ 38 (2001), 37--100; II. {\em On
weighted projective planes and their affine rulings}. {\it Ibid.},
101--150; III. {\em On log $\Bbb Q$-homology planes and weighted
projective planes}, Canad.\ J.\ Math.\ 56 (2004), 1145--1189.

\bibitem[De]{De} A.\ S.\ Demushkin,
{\em Combinatorial invariance of toric singularities}. Vestnik
Moskov.\ Univ.\ Ser. I Mat.\ Mekh.\ 2 (1982), 80--87.

\bibitem[Du$_1$]{Du1} A.\ Dubouloz,
{\em Completions of normal affine surfaces with a trivial
Makar-Limanov invariant}. Michigan Math.\ J.\ 52 (2004), 289--308.

\bibitem[Du$_2$]{Du2} A.\ Dubouloz, I.
{\em Generalized Danielewski surfaces}. Preprint, Institut Fourier
des math\'ematiques {\bf 612} (2003), math.AG/0401225 (2004);
Transfor. Groups 10 (2005), 139-162. II. {\em Embeddings of
generalized Danielewski surfaces in affine spaces}.
math.AG/0403208 (2004); Compos.\ Math.\ ({\it to appear}).

\bibitem[Du$_3$]{Du3} A.\ Dubouloz, {\em Sur une classe de sch\'emas
avec actions de fibr\'es en droites}. PhD thesis, Grenoble, 2004,
http://tel.ccsd.cnrs.fr/documents/archives0/00/00/77/33/index.html.

\bibitem[FiKa]{FiKa} K.-H.\ Fieseler, L.\ Kaup, {\em
Hyperbolic $\C\sp *$-actions on affine algebraic surfaces}.
Complex analysis (Wuppertal, 1991), 160-168, Aspects Math.\ E17,
Vieweg, Braunschweig, 1991.

\bibitem[FKZ]{FKZ}
H.\ Flenner, S.\ Kaliman, M.\ Zaidenberg, {\em Birational
transformations of weighted graphs} ({\it in this volume}).


\bibitem[FlZa$_1$]{FlZa1}
H.\ Flenner, M.\ Zaidenberg, I. {\em Normal affine surfaces with
$\C^*$-actions}. Osaka J.\ Math.\ 40, 2003,  981--1009.
%
II. {\em Locally nilpotent derivations on affine surfaces with a
$\C^*$-action}. math.AG/0403215; Osaka J.\ Math.\ 42:4, 2005 ({\it
to appear}).

\bibitem[FlZa$_2$]{FlZa2} H.\ Flenner, M.\ Zaidenberg,
{\em On the uniqueness of $\C^*$-actions on affine surfaces}.
"Affine Algebraic Geometry" (conference Proceedings volume, Jaime
Gutierrez, Vladimir Shpilrain, and Jie-Tai Yu, Eds.),
Contemporary Mathematics 369, Amer.\ Math.\ Soc.\, Providence,
R.I.\ 2005, 97-111.

\bibitem[Fu]{Fu} T. Fujita, {\em On the topology of noncomplete
algebraic surfaces}. J.\ Fac.\ Sci.\ Univ.\ Tokyo Sect.\ IA Math.\
29 (1982), 503--566.

\bibitem[Gi]{Gi} M.\ H.\ Gizatullin, I.
{\em Affine surfaces that are quasihomogeneous with respect to an
algebraic group}. Math. USSR Izv. 5 (1971), 754--769; II. {\em
Quasihomogeneous affine surfaces}. {\it Ibid.} 1057--1081.

\bibitem[Gr]{Gr}
H.\ Grauert, {\em \"Uber Modifikationen und exzeptionelle
analytische Mengen}. Math.\ Ann.\ {\em 146} (1962), 331--368.

\bibitem[Gu]{Gu} J.\ Gubeladze, {\em The isomorphism problem
for commutative
monoid rings}. J.\ Pure Appl.\ Algebra 129 (1998), 35--65.

\bibitem[GMMR]{GMMR} R.\ V.\ Gurjuar, K.\ Masuda, M.\
Miyanishi and P.\ Russell, {\em Affine lines on affine surfaces
and the Makar-Limanov invariant}. Preprint 2005, 42p.


\bibitem[Hi$_1$]{Hi1} F.\ Hirzebruch,
{\em \"Uber vierdimensionale Riemannsche Flächen mehrdeutiger
analytischer Funktionen von zwei komplexen Veränderlichen}. Math.\
Ann.\ 126 (1953), 1--22.

\bibitem[Hi$_2$]{Hi2} F.\ Hirzebruch,
{\em The topology of normal singularities of an algebraic surface
(after D. Mumford)}. S\'eminaire Bourbaki, Vol. 8, Exp. No. 250,
129--137, Soc. Math. France, Paris, 1995.


\bibitem[ML$_1$]{ML1} L.\ Makar-Limanov, {\em On groups of
automorphisms
of a class of surfaces}. Israel J.\ Math.\ 69 (1990), 250--256;
II. {\em  On the group of automorphisms of a surface $x^ny=P(z)$},
Israel J.\ Math.\ 121 (2001), 113--123.

\bibitem[ML$_2$]{ML2} L.\ Makar-Limanov,
{\em Locally nilpotent derivations on the surface $xy=p(z)$}.
Proceedings of the Third International Algebra Conference
(Tainan, 2002), Kluwer Acad.\ Publ.\, Dordrecht, 2003, 215--219.

\bibitem[MaMi$_1$]{MaMi1}
K.\ Masuda, M.\ Miyanishi, {\em  The additive group actions on
$\Q$-homology planes}. Annales de l'Institut Fourier  53 (2003),
429-464.

\bibitem[MaMi$_2$]{MaMi2}
K.\ Masuda, M.\ Miyanishi, {\em  Affine pseudo-planes with torus
actions}. Preprint, 2005, 22p.

\bibitem[Mi]{Mi} M.\ Miyanishi, {\em Open algebraic surfaces}.
CRM Monograph Series, 12. American Mathematical Society,
Providence, RI, 2001.

\bibitem[Mu]{Mu}
D.\ Mumford, {\em The topology of normal singularities of an
algebraic surface and a criterion for simplicity}. Inst.\ Hautes
\'Etudes Sci.\ Publ.\ Math.\ 9 (1961), 5--22.



\bibitem[OrWa]{OrWa} P.\ Orlik, P.\ Wagreich, {\em Equivariant
resolution
of singularities with $C\sp{*}$ action}. Proceedings of the
Second Conference on Compact Transformation Groups (Univ.\
Massachusetts, Amherst, Mass., 1971), Part I, 270--290. Lecture
Notes in Math.\ 298, Springer, Berlin, 1972.


\bibitem[Ru]{Rus} P.\ Russell, {\em Some formal aspects of
the theorems of Mumford-Ramanujam}. Algebra, arithmetic and
geometry, Part I, II (Mumbai, 2000), Tata Inst.\ Fund.\ Res.\
Stud.\ Math.\, 16, Tata Inst.\ Fund.\ Res.\, Bombay, 2002,
557--584.

\bibitem[Sa]{Sa}  F.\ Sakai, {\em Weil divisors on normal
surfaces}. Duke Math.\ J.\ 51 (1984), 877--887.

\bibitem[Sch]{Sch} S. Schr\"oer, I. {\em On non-projective normal
surfaces}. Manuscripta Math. 100 (1999), 317--321. II. {\em On
contractible curves on normal surfaces}. J.\ Reine Angew.\ Math.\
524 (2000), 1--15.

\bibitem[Su]{Su} H.\ Sumihiro, {\em  Equivariant completion}.
I. J.\ Math.\ Kyoto Univ.\ 14 (1974), 1--28; II. {\it ibid.} 15
(1975), 573--605.

\bibitem[Ta]{Ta} R.\ Tanimoto, {\em Algebraic torus actions
on affine algebraic surfaces}. J.\ Algebra 285 (2005), 73--97.

\bibitem[Zar]{Zar} O.\ Zariski, {\em The reduction of the
singularities of an algebraic surface}. Ann.\ Math.\ (2) 40,
639-689 (1939).


\end{thebibliography}
\end{document}